\documentclass[10pt,draftcls,onecolumn]{IEEEtran}

\IEEEoverridecommandlockouts                                 

\newtheorem{theorem}{Theorem}
\newtheorem{corollary}[theorem]{Corollary}
\newtheorem{lemma}[theorem]{Lemma}
\newtheorem{proposition}[theorem]{Proposition}

\newtheorem{remark}{Remark}

\newcommand{\enma}[1]   {\ensuremath{#1}}

\newcommand{\non}{\nonumber}

\newcommand{\beq}{\begin{equation}}
\newcommand{\eeq}{\end{equation}}
\newcommand{\bseq}{\begin{subequations}}
\newcommand{\eseq}{\end{subequations}}
\newcommand{\beqn}{\begin{eqnarray}}
\newcommand{\eeqn}{\end{eqnarray}}
\newcommand{\ba}{\begin{array}}
\newcommand{\ea}{\end{array}}
\newcommand{\bct}{\begin{center}}
\newcommand{\ect}{\end{center}}
\newcommand{\btmz}{\begin{itemize}}
\newcommand{\etmz}{\end{itemize}}
\newcommand{\benum}{\begin{enumerate}}
\newcommand{\eenum}{\end{enumerate}}

\newcommand{\cH}{\enma{\mathcal H}}
\newcommand{\cL}{\enma{\mathcal L}}

\newcommand{\norm}[1]{\| #1 \|}                 

\newcommand{\diag}      {\enma{\mathrm{diag}}}

\newcommand{\inner}[2]{\left\langle #1,#2 \right\rangle}

\newcommand{\matbegin}{
        \left[
}
\newcommand{\matend}{
        \right]
}

\newcommand{\tbo}[2]{
  \matbegin \begin{array}{c}
       #1 \\ #2
       \end{array} \matend }

\newcommand{\obt}[2]{
  \matbegin \begin{array}{cc}
       #1 & #2
       \end{array} \matend }
\newcommand{\obth}[3]{
  \matbegin \begin{array}{ccc}
       #1 & #2 & #3
       \end{array} \matend }
\newcommand{\tbt}[4]{
  \matbegin \begin{array}{cc}
       #1 & #2 \\ #3 & #4
       \end{array} \matend }

\newcommand{\thbth}[9]{
 \matbegin \begin{array}{ccc}
                #1 & #2 & #3 \\
                #4 & #5 & #6 \\
                #7 & #8 & #9
                \end{array}\matend}

\newcommand{\fbo}[4]{
 \matbegin \begin{array}{c}
                #1 \\ #2 \\ #3 \\ #4
                \end{array}\matend}

\newcommand{\be}{\begin{equation}}
\newcommand{\ee}{\end{equation}}

\newcommand{\cplxs}{ C\kern -.35em \rule{0.03 em}{.7 ex}~   }

\def\complex{\hbox{C\kern -.45em \rule{0.03 em}{1.5 ex}}~}

\newcommand{\bi}{\begin{itemize}}
\newcommand{\ei}{\end{itemize}}

\newtheorem{assumption}{Assumption}

\usepackage{amssymb,latexsym,color,amsmath,epsfig,graphicx,amsmath}

\usepackage[usenames,dvipsnames,svgnames,table]{xcolor}
\usepackage{tikz}
\usepackage{subfig}
\usepackage{mhchem}
\usetikzlibrary{calc,arrows,automata,positioning,shapes}
\usepackage{pgfplots}
\pgfplotsset{compat=newest}
\pgfplotsset{plot coordinates/math parser=false}
\usepgfplotslibrary{patchplots}

\usepackage{setspace,cite,dsfont} \usepackage{cite}
\usepackage{algorithm,algorithmic}
\usepackage{hyperref}
\usepackage{listings}

\usepackage{accents}
\newcommand{\cC}{{\cal C}}

\newcommand{\cM}{{\cal V}}

\newcommand{\cV}{{\cal V}}
\newcommand{\cS}{{\cal S}}

\newcommand{\mre}{\mathrm{e}}
\newcommand{\bbR}{\mathbb{R}}

\newcommand{\ds}{\displaystyle}

\newcommand{\tc}{\textcolor}

\definecolor{matlabblue}{rgb}{0   0.4470   0.7410}
\definecolor{matlabred}{rgb}{0.8500    0.3250    0.0980}
\definecolor{matlaborange}{rgb}{0.9290    0.6940    0.1250}
\newcommand{\vsp}{\vspace*{0.1cm}}

\DeclareMathOperator*{\argmin}{argmin}
\DeclareMathOperator*{\minimize}{minimize}
\DeclareMathOperator*{\subject}{subject~to}

\newcommand{\mrd}{\mathrm{d}}

\newcommand{\mrj}{\mathrm{j}}

\newcommand{\zms}{z_\mu^\star}

\newcommand{\DefinedAs}[0]{\mathrel{\mathop:}=}

\newcommand{\sign}{\mathrm{sign}}

\newcommand{\tx}{{\tilde{x}}}
\newcommand{\ty}{{\tilde{y}}}
\newcommand{\tw}{{\tilde{w}}}

\newcommand{\prox}{\mathbf{prox}}

\newcommand{\fx}{f}
\newcommand{\gz}{g}

\newcommand{\mor}{M}
\newcommand{\gr}{\nabla f}
\newcommand{\Hes}{H}
 
 \newcommand{\tcbl}{\tc{black}}
 
 \newcommand{\bbP}{\mathbb{P}}
\DeclareMathOperator*{\blkdiag}{blkdiag}

\begin{document}

\title{A second order primal-dual method for \\[-0.35cm] nonsmooth convex composite optimization}
\author{Neil K.\ Dhingra, Sei Zhen Khong, and Mihailo R.\ Jovanovi\'c
\thanks{Financial support from the National Science Foundation under awards ECCS-1708906 and ECCS-1809833 is gratefully acknowledged.}
\thanks{N.\ K.\ Dhingra is with Numerica Corporation, Fort Collins, CO 80528. S.\ Z.\ Khong is an independent researcher. M.\ R.\ Jovanovi\'c is with the Ming Hsieh Department of Electrical and Computer Engineering, University of Southern California, Los Angeles, CA 90089. E-mails: neil.k.dh@gmail.com, szkhongwork@gmail.com, mihailo@usc.edu.}}

\maketitle

	\vspace*{-0.75cm}

	\begin{abstract}
We develop a second order primal-dual method for optimization problems in which the objective function is given by the sum of a strongly convex twice differentiable term and a possibly nondifferentiable convex regularizer. After introducing an auxiliary variable, we utilize the proximal operator of the nonsmooth regularizer to transform the associated augmented Lagrangian into a function that is once, but not twice, continuously differentiable. The saddle point of this function corresponds to the solution of the original optimization problem. We employ a generalization of the Hessian to define second order updates on this function and prove global exponential stability of the corresponding differential inclusion. Furthermore, we develop a globally convergent customized algorithm that utilizes the primal-dual augmented Lagrangian as a merit function. We show that the search direction can be computed efficiently and prove quadratic/superlinear asymptotic convergence. We use the $\ell_1$-regularized \tcbl{model predictive control} problem and the problem of designing a distributed controller for a spatially-invariant system to demonstrate the merits and the effectiveness of our method.
	\end{abstract}

	\vspace*{-3ex}
\section{Introduction}

We study a class of composite optimization problems in which the objective function is given by the sum of a differentiable strongly convex component and a nondifferentiable convex component. Problems of this form are encountered in diverse fields including 
compressive sensing,
 machine learning, statistics, image processing, and control~\cite{don06,golosh09,charecparwil12,biclev08,linfarjovTAC13admm,jovdhiEJC16,dhicoljovTCNS19,mogdhijovCDC16,galmac12}. They often arise in structured feedback synthesis where it is desired to balance controller performance (e.g., the closed-loop ${\cal H}_2$ or ${\cal H}_\infty$ norm) with structural complexity~\cite{linfarjovTAC13admm,jovdhiEJC16}.

	The lack of differentiability in the regularization term precludes the use of standard descent methods for smooth objective functions. Proximal gradient methods~\cite{liomer79,compes11,parboy13} and their accelerated variants~\cite{bectebSIAM09} generalize gradient descent, but typically require the nonsmooth term to be separable.
	
	An alternative approach introduces an auxiliary variable to split the smooth and nonsmooth components of the objective function. The reformulated problem facilitates the use of splitting methods such as the alternating direction method of multipliers (ADMM)~\cite{boyparchupeleck11}. This augmented-Lagrangian-based method divides the optimization problem into simpler subproblems, allows for a broader class of regularizers than proximal gradient, and it is convenient for distributed implementation. In~\cite{dhikhojovTAC19}, we exploited the structure of proximal operators associated with nonsmooth regularizers and introduced the {\em proximal augmented Lagrangian\/}. This continuously differentiable function enables the use of the standard method of multipliers (MM) for nonsmooth optimization and it is convenient for distributed implementation via the Arrow-Hurwicz-Uzawa gradient flow dynamics. Recent work has extended MM to incorporate second order updates of the primal and dual variables~\cite{gilrobCOA12,armbenomhpat14,armomh15} for nonconvex problems with twice continuously differentiable objective functions.
		
Since first order approaches tend to converge slowly to high-accuracy solutions, much work has focused on developing second order methods for nonsmooth composite optimization. A generalization of Newton's method was developed in~\cite{becfad12,tseyu09,byrnocozt15,leesunsau14} and it requires solving a regularized quadratic subproblem to determine a search direction. Related ideas have been utilized for sparse inverse covariance estimation in graphical models~\cite{hsidhiravsus11,hsisusdhirav14} and topology design in \mbox{consensus networks~\cite{mogjovTCNS18}.}
 
Generalized Newton updates for identifying stationary points of (strongly) semismooth gradient mappings were first considered in~\cite{qixsun93,qixsun99,mif77}. In~\cite{patstebem14,stethepat17,thestepat18}, the authors introduce the once-continuously differentiable Forward-Backward Envelope (FBE) and solve composite problems by minimizing the FBE using line search, quasi-Newton methods, or second order updates based on an approximation of the generalized Hessian.

We develop a second order primal-dual algorithm for nonsmooth composite optimization by leveraging these advances. Second order updates for the once continuously differentiable proximal augmented Lagrangian are formed using a generalized Hessian. We employ the merit function introduced in~\cite{gilrobCOA12} to assess progress towards the optimal solution and develop a globally convergent customized algorithm with fast asymptotic convergence rate. When the proximal operator associated with a nonsmooth regularizer is (strongly) semismooth, our algorithm exhibits local (quadratic) superlinear convergence. 

Our presentation is organized as follows. In Section~\ref{sec.problem}, we formulate the problem and provide necessary background material. In Section~\ref{sec.2nd}, we define second order updates to find the saddle points of the proximal augmented Lagrangian. In Section~\ref{sec.diffinc}, we prove global exponential stability of a continuous-time differential inclusion. In Section~\ref{sec.pd}, we develop a customized algorithm that converges globally to the saddle point and exhibits superlinear or quadratic asymptotic convergence rates. In Section~\ref{sec.ex}, we provide examples to illustrate the utility of our method. We discuss interpretations of our approach and connections to the alternative algorithms in Section~\ref{sec.disc} and conclude the paper in Section~\ref{sec.conc}.

	\vspace*{-2ex}
\section{Problem formulation and background}
	\label{sec.problem}

We consider the problem of minimizing the sum of two convex functions over an optimization variable $x \in \bbR^n$,
\be
	\label{pr}
	\minimize\limits_x ~~ \fx(x) \;+\; \gz(Tx)
\ee
where $T \in \bbR^{m \times n}$. Problem~\eqref{pr} \tcbl{is widely used} 
in compressive sensing to incorporate structural considerations into traditional sensing or regression~\cite{don06,golosh09,charecparwil12}. \tcbl{It has also gained traction in many fields including control-theoretic applications~\cite{jovdhiEJC16,galmac12,dhicoljovTCNS19,linfarjovTAC13admm,mogdhijovCDC16}.}

\begin{assumption} 
\label{ass.func}
The function $f$ is twice continuously differentiable, has an $L_f$ Lipschitz continuous gradient $\gr$, and is strictly convex with $\nabla^2 f \succ 0$; the function $\gz$ is proper, \tcbl{closed}, and convex; and the matrix $T$ has full row rank. 
\end{assumption}

As specified in Assumption~\ref{ass.func}, the differentiable part of the objective function $f$, which quantifies performance loss, is strictly convex with a Lipschitz continuous gradient. In contrast, the regularization function $g$ may be nondifferentiable and is used to incorporate structural requirements on the optimization variable. In structured feedback synthesis~\cite{linfarjovTAC13admm,jovdhiEJC16}, $f$ typically quantifies the closed-loop performance, e.g., the $\cH_2$ norm, and $g$ imposes structural requirements  on the controller, e.g., by penalizing the amount of network traffic~\cite{dorjovchebulACC13,dorjovchebulTPS14,wudorjovTPS16,wujovSCL17}. Although our strongest results require strong convexity of $f$, our theory and techniques are applicable as long as the Hessian of $f$ is positive definite. \tcbl{The extension to weakly convex functions using the pseudoinverse of the Hessian or assuming ``restricted strong convexity'' is currently under study.}

The matrix $T$ is important when the desired structure has a simple representation in the co-domain of $T$, but it makes the problem more challenging. One approach is to reformulate~\eqref{pr} by introducing an auxiliary optimization variable $z \in \bbR^m$,
	\be
	\label{pr.split}
	\ba{rl}
	\minimize\limits_{x, \, z} & \fx(x) \;+\; \gz(z) 
		\\[0.1cm]
	\subject    & Tx \;-\; z \;=\; 0.
	\ea
	\ee
Problem~\eqref{pr.split} is convenient for constrained optimization algorithms based on the augmented Lagrangian,
\be
	\label{eq.al}
	\non
	\cL_\mu(x,z;y)
	\DefinedAs
	\fx(x) + \gz(z) + y^T(Tx -z) + \tfrac{1}{2\mu} \, \norm{Tx-z}^2,
\ee
where $y \in \bbR^m$ is the Lagrange multiplier and $\mu$ is a positive parameter. Relative to the standard Lagrangian, $\cL_\mu$ contains an additional quadratic penalty on the linear constraint in~\eqref{pr.split}.

In the remainder of this section, we provide background on proximal operators and describe generalizations of the gradient for nondifferentiable functions. We also briefly overview existing approaches for solving~\eqref{pr}. 

	\vspace*{-3ex}
\subsection{Background}

\subsubsection{Proximal operators}
	\label{sec.back_back}

The proximal operator of the function $g$ is the minimizer of the sum of $g$ and a proximal term,
\begin{subequations}
\be
	\prox_{\mu \gz}(v)
	\;\DefinedAs\;
	\argmin_z 
	\,
	\left(
	\gz(z)
	\;+\;
	\tfrac{1}{2\mu} \, \norm{z \,-\, v}^2
	\right)
		\label{eq.prox}
\ee
where $\mu$ is a positive parameter and $v$ is a given vector. When $g$ is convex, its proximal operator is Lipschitz continuous with parameter $1$, differentiable almost everywhere, and firmly non-expansive~\cite{parboy13}. The value function associated with~\eqref{eq.prox} specifies the Moreau envelope of $g$,
	\beq
	\label{eq.mor}
	\!\!\!\!
	\ba{rrl} 
	M_{\mu g} (v) 
	&\!\!\! =\!\!\! &
	g (\prox_{\mu g}(v)) 
	\, + \, 
	\tfrac{1}{2\mu} \, \norm{\prox_{\mu g}(v) - v}^2.
	\ea
	\eeq
The Moreau envelope is continuously differentiable, even when $g$ is not, and its gradient
\be
	\nabla M_{\mu g}(v)
	\;=\;
	\tfrac{1}{\mu}
	\left(v \;-\; \prox_{\mu g}(v) \right)
	\label{eq.gradM}
\ee
\end{subequations}
is Lipschitz continuous with parameter $1/\mu$. 

For example, the proximal operator associated with the $\ell_1$ norm, $g (z) = \sum |z_i|$, is determined by soft-thresholding,
	$
	\cS_{\mu} (v_i)
	\DefinedAs 
         \sign (v_i)
         \max\{|v_i| - \mu, 0\},
         $
the associated Moreau envelope is the Huber function,
	$
    M_{\mu g}(v_i)
    = 
    \{
    v_i^2 / (2 \mu),
    \,
    |v_i|
    \leq
    \mu;
    $
    $
    |v_i| - \mu/2,
    \,
    |v_i|
    \geq
    \mu
    \},
    $
and its gradient is the saturation function,
	$
    \nabla M_{\mu g}(v_i)
    =
    \sign(v_i) \min \{ |v_i|/\mu, \,1 \}.
    $
Other regularizers with efficiently computable proximal operators include indicator functions of simple convex sets as well as the nuclear and Frobenius norms of matricial variables. Such regularizers can enforce bounds on $x$, promote low rank solutions, and enhance group sparsity, respectively. 
 
\subsubsection{Generalization of the gradient and Jacobian}
Although $\prox_{\mu g}$ is typically not differentiable, it is Lipschitz continuous and therefore differentiable almost everywhere~\cite{rocwet09}. One generalization of the gradient for such functions is given by the $B$-subdifferential set~\cite{rob87}, which applies to locally Lipschitz continuous functions $h$: $\bbR^m \to \bbR$. Let $C_h$ be a set at which $h$ is differentiable. Each element in the set $\partial_B h (\bar{z})$ is the limit point of a sequence of gradients $\{\nabla h(z_k)\}$ evaluated at a sequence of points $\{z_k\} \subset C_h$ whose limit is $\bar z$,
\be
	\label{eq.Bsub}
	\partial_B h(\bar z)
	\; \DefinedAs \, 
	\left\{
	J
	\,|~
	\exists \{z_k\}
	\subset
	C_h,
	~
	z_k \to \bar z,
	~
	\nabla h(z_k) \to J
	\right\}.
\ee
{If $h$ is continuously differentiable in the neighborhood of a point $z$, the $B$-subdifferential set becomes single valued and it is given by the gradient, $\partial_B h (z) = \nabla h (z)$. In general, $\partial_B h(\bar z)$ is not a convex set; e.g., if $h (z) = | z |$, $\partial_B h(0) = \{ -1, 1 \}$.}

The Clarke subdifferential set of $h$: $\bbR^m \to \bbR$ at $\bar z$ is the convex hull of the $B$-subdifferential set~\cite{cla90},
	$
	\partial_C h(\bar z)
	\DefinedAs
	\mbox{conv} 
	\,
	(
	\partial_B h(\bar z)
	).
	$
When $h$ is a convex function, the Clarke subdifferential set is equal to the subdifferential set $\partial h(\bar z)$ which defines the supporting hyperplanes of $h$ at $\bar z$~\cite[Chapter VI]{hrilem13}. For a function $G$: $\bbR^m \to \bbR^n$, the $B$-generalization of the Jacobian at a point $\bar z$ is given by
\[
	\partial_B G(\bar z)
	\;\DefinedAs\;
	\obth{J_1^T}{\dots}{J_n^T}^T
\]
where each $J_i \in \partial_B G_i(\bar z)$ is a member of the $B$-subdifferential set of the $i$th component of $G$ evaluated at $\bar z$. The Clarke generalization of the Jacobian at a point $\bar z$, $\partial_CG(\bar z)$, has the same structure where each $J_i \in \partial_C G_i(\bar z)$ is a member of the Clarke subdifferential of $G_i(\bar z)$.

\subsubsection{Semismoothness}

The mapping $G$: $\bbR^m \to \bbR^n$ is semismooth at $\bar z$ if for any sequence $z_k \to \bar z$, the sequence of Clarke generalized Jacobians $J_{G_k} \in \partial_C G(z_k)$ provides a first order approximation of $G$,
\be
	\norm{G(z_k) \,-\, G(\bar z) \,+\, J_{G_k} (\bar z \,-\, z_k)}
	\;=\;
	o(\norm{z_k \,-\, \bar z})
	\label{eq.semi}
\ee
where $\phi(k) = o(\psi(k))$ denotes that $\phi(k)/\psi(k) \to 0$ as $k$ tends to infinity~\cite{harwri79}.
The function $G$ is strongly semismooth if this approximation satisfies the stronger condition,
\[
	\norm{G(z_k) \,-\, G(\bar z) \,+\, J_{G_k} (\bar z \,-\, z_k)}
	\;=\;
	O(\norm{z_k \,-\, \bar z}^2),
\]
where $\phi(k) = O(\psi(k))$ signifies that $|\phi(k)| \leq L \psi(k)$ for some positive constant $L$ and positive $\psi(k)$~\cite{harwri79}.

	\vsp
	
\begin{remark}
	\label{sec.semi}
	(Strong) semismoothness of the proximal operator leads to fast asymptotic convergence of the differential inclusion (see Section~\ref{sec.diffinc}) and the efficient algorithm (see Section~\ref{sec.pd}). Proximal operators associated with many typical regularization functions (e.g., the $\ell_1$ and nuclear norms~\cite{jiasuntoh14}, piecewise quadratic functions~\cite{mensunzha05}, and indicator functions of affine convex sets~\cite{mensunzha05}) are strongly semismooth. In general, semismoothness of $\prox_{\mu g}$ follows from semismoothness of the projection onto the epigraph of $g$~\cite{mensunzha05}. However, there are convex sets onto which projection is not directionally differentiable~\cite{kru69}. The indicator functions associated with such sets or functions whose epigraph is described by such sets may induce proximal operators which are not semismooth.
	\end{remark}

	\vspace*{-3ex}
\subsection{Existing methods}

Problem~\eqref{pr} is encountered in a host of applications and it has been the subject of extensive study. Herein, we provide a brief overview of existing approaches to solving it. 

	\vsp
	
\subsubsection{First order methods}
When $T$ is identity or a diagonal matrix, the proximal gradient method, which generalizes gradient descent to certain classes of nonsmooth composite optimization problems~\cite{bectebSIAM09,parboy13}, can be used to solve~\eqref{pr},
\[
	x^{k+1}
	\;=\;
	\prox_{\alpha^k \gz}
	\!
	\left(x^k \,-\, \alpha^k \gr(x^k)\right)
\]
where $x^k$ is the current iterate and $\alpha^k$ is the step size. When $\gz = 0$, we recover gradient descent, when $\gz$ is the indicator function $I_\cC(x)$ of the convex set $\cC$, it simplifies to projected gradient descent, and when $g$ is the $\ell_1$ norm, it corresponds to the Iterative Soft-Thresholding Algorithm (ISTA). Nesterov-style techniques can also be employed for acceleration~\cite{bectebSIAM09}. \tcbl{When $T$ is not a diagonal matrix, proximal methods are harder to use because efficient computation of the proximal operator associated with $g \circ T$ may not follow from efficient computation of $\prox_{\mu \gz}$. However, if the conjugate function of $f^*(y) \DefinedAs \max_x\{\inner{y}{x} - f(x)\}$, its gradient, and the proximal operator associated with the conjugate function of $g$ are known or easily computed, proximal methods on the dual problem can be employed to solve~\eqref{pr}~\cite{becteb14,tse91,gisboy15}; these approaches are particularly useful when the function $f$ is separable.}

The alternating direction method of multipliers (ADMM) \tcbl{can be used to solve~\eqref{pr} via~\eqref{pr.split} even when the matrix $T$ is not diagonal.} \tcbl{ADMM entails} alternating between minimization of $\cL_\mu(x,z;y)$ over $x$ (a continuously differentiable problem), minimization of $\cL_\mu(x,z;y)$ over $z$ (amounts to evaluating $\prox_{\mu g}$), and a gradient ascent step in $y$~\cite{boyparchupeleck11},
\be
	\ba{rcl}
	x^{k+1}
	&\!\!=\!\!&
	\ds\argmin_x \cL_\mu(x,z^k;y^k)
	\\[0.15cm]
	z^{k+1}
	&\!\!=\!\!&
	\prox_{\mu \gz} (Tx^{k+1} \,+\, \mu y^k)
	\\[0.15cm]
	y^{k+1}
	&\!\!=\!\!&
	y^k \,+\, \tfrac{1}{\mu} \, (Tx^{k+1} \,-\, z^{k+1}).
	\ea
	\label{eq.ADMM}
\ee
\tcbl{A similar method is the alternating minimization algorithm (AMA), whose form is identical to the ADMM iteration~\eqref{eq.ADMM} except that the first step requires minimization of the {\em standard\/} Lagrangian over $x$~\cite{tse91}. We note that although second order methods may be used to solve {\em inner\/} $x$-minimization problem, the overall {\em outer\/} iteration in~\eqref{eq.ADMM} corresponds to a sequence of first-order updates of the dual variable $y$.}

	\vsp
	
\subsubsection{Second order methods}

The slow convergence of first order methods to high-accuracy solutions motivates the development of second order methods for solving~\eqref{pr}. A generalization of Newton's method for nonsmooth problems~\eqref{pr} with $T = I$ was developed in~\cite{tseyu09,becfad12,leesunsau14,byrnocozt15}. A sequential quadratic approximation of the smooth part of the objective function is utilized and a search direction $\tilde{x}$ is obtained as the solution of a regularized quadratic subproblem,
\be
	\minimize_\tx
	~~
	\tfrac{1}{2}
	\,
	\tx^TH \tx
	\;+\;
	\nabla f(x^k)^T\tx
	\;+\;
	g(x^k + \, \tx)
	\label{eq.proxnewt}
\ee
where $x^k$ is the current iterate and $H$ is the Hessian of $f$. This method generalizes the projected Newton method~\cite{ber82pn} to a broader class of regularizers. For example, when $g$ is the $\ell_1$ norm, this amounts to solving a LASSO problem~\cite{tib96}, which can be a challenging task. Coordinate descent is often used to solve this subproblem~\cite{leesunsau14} and it has been observed to perform well in practice~\cite{hsidhiravsus11,hsisusdhirav14,mogjovTCNS18}.

The Forward-Backward Envelope (FBE) was introduced in~\cite{patstebem14,stethepat17,thestepat18}. The FBE is a once-continuously differentiable nonconvex function of $x$ and its minimum corresponds to the solution of~\eqref{pr} with $T = I$. As demonstrated in Section~\ref{sec.disc}, FBE can be obtained from the proximal augmented Lagrangian (that we introduce in Section~\ref{sec.2nd}). Since the generalized Hessian of FBE involves third-order derivatives of $f$ (which may be expensive to compute), references~\cite{patstebem14,stethepat17,thestepat18} employ either truncated- or quasi-Newton methods to obtain a second order update to $x$. \tcbl{Generalizations of the FBE to the problems with $T \neq I$ have also been developed~\cite{stethepat18}.}

	\vspace*{-2ex}
\section{The proximal augmented Lagrangian and second order updates}
	\label{sec.2nd}

	In this section, we transform $\cL_\mu(x,z;y)$ into a form that is once but not twice continuously differentiable. \tcbl{Under Assumption~\ref{ass.func},} the resulting function, which we call the proximal augmented Lagrangian, \tcbl{is convex in $x$ and concave in $y$ and its saddle points correspond to the solution of~\eqref{pr}. Even though methods of this paper can be also employed for nonconvex problems, in this case, a saddle point of the proximal augmented Lagrangian identifies a local minimum of~\eqref{pr} and convergence is not guaranteed. We} define second order updates to find its saddle points, show that they are always well defined, and prove that they are locally (quadratically) superlinearly convergent when $\prox_{\mu g}$ is (strongly) semismooth.

	\vspace*{-2ex}
\subsection{Proximal augmented Lagrangian}
	\label{sec.pal}

The continuously differentiable proximal augmented Lagrangian was recently introduced in~\cite{dhikhojovTAC19}. This was done by rewriting $\cL_\mu (x,z;y)$ via completion of squares,
    \[
    \cL_\mu(x,z;y)
    \, = \,
    f(x)
    \, + \,
    g(z)
    \, + \,
    \tfrac{1}{2\mu} \, \norm{z \, - \, (Tx + \mu y)}^2
    \, - \,
    \tfrac{\mu}{2} \, \norm{y}^2
    \] 
and restricting it to the manifold that corresponds to explicit minimization over the auxiliary variable $z$,
	\be
	\cL_\mu (x; y)
	\; \DefinedAs \;
	\cL_\mu(x, \zms(x,y); y)
	\non
	\ee
where the minimizer of $\cL_\mu(x,z;y)$ over $z$ is determined by the proximal operator of the function $g$,
	\be
	\zms(x,y)
	\; = \;
	\argmin\limits_z
	 \,
	 \cL_\mu(x, z; y)
	\; = \;
	\prox_{\mu g}(Tx \, + \, \mu y)
	\non
	\ee
and it defines the aforementioned manifold.
	
	\vsp

\begin{theorem}[Theorem 1 in~\cite{dhikhojovTAC19}] \label{thm.diff}
Let Assumption~\ref{ass.func} hold. Then, minimization of the augmented Lagrangian $\cL_\mu(x,z; y)$ associated with problem~\eqref{pr.split} over $(x,z)$ is equivalent to minimization of the {\em proximal augmented Lagrangian\/}
\be
	\cL_\mu(x;y)
	\;=\;
	f(x)
    \; + \;
    M_{\mu g}
    (Tx \, + \, \mu y)
    \; - \;
    \tfrac{\mu}{2} \, \| y \|^2
    \label{eq.alprox}
\ee
over $x$. Moreover, $\cL_\mu(x;y)$ is continuously differentiable over $x$ and $y$ and its gradient,
	\be
	\nabla \cL_\mu (x;y)
	\; = \;
	\tbo{\nabla_x \cL_\mu (x;y)}{\nabla_y \cL_\mu (x;y)}
	\; = \;
	\tbo{\gr(x) \, + \, T^T\nabla M_{\mu g}(T x  + \mu y)}
	{\mu\nabla M_{\mu g}(T x  + \mu y) \,-\, \mu y}
	\label{eq.gradLmu}
\ee
is Lipschitz continuous.
\end{theorem}

	\vsp

The proximal augmented Lagrangian $\cL_\mu(x;y)$ contains the Moreau envelope of $g$ and its introduction allows the use of the method of multipliers (MM) to solve problem~\eqref{pr.split}. MM requires {\em joint\/} minimization of $\cL_\mu(x,z;y)$ over $x$ and $z$ which is, in general, challenging because the $(x,z)$-minimization subproblem is nondifferentiable. However, Theorem~\ref{thm.diff} enables an equivalent implementation of MM
	\be
	\label{eq.MM}
	\ba{rcl}
	x^{k+1}
	& \!\! = \!\! &
	\argmin\limits_{x}
	\,
	\cL_\mu (x; y^k)
	\\[0.3cm]	
	y^{k+1}
	& \!\! = \!\! &
	y^k
	~+~
	\tfrac{1}{\mu} \,
	(Tx^{k+1} \, - \, z^\star_\mu (x^{k+1},y^k) )
	\ea
	\ee
which improves performance relative to ADMM~\cite{dhikhojovTAC19}. 
	
Continuous differentiability of $\cL_\mu(x;y)$ also enables joint \tcbl{gradient} updates of $x$ and $y$ via the primal-dual Arrow-Hurwicz-Uzawa iteration,
\be
	\ba{rcl}
    \tcbl{x^{k+1}}
	&\!\!=\!\!&
    \tcbl{x^k
        \; - \; 
        \alpha_k}
	\nabla_x \cL_\mu(x;y)
	\\[0.15cm]
    \tcbl{y^{k+1}}
	&\!\!=\!\!&
    \tcbl{y^k
        \; + \;
        \alpha_k}
	\nabla_y \cL_\mu(x;y)
	\ea
    \label{eq.ahu}
\ee
where $\nabla_x \cL_\mu$ and $\nabla_y \cL_\mu$ are given by~\eqref{eq.gradLmu} \tcbl{and $\alpha_k$ is stepsize}. When $\nabla f$ and $T$ are sparse mappings, this method is convenient for distributed implementation and it is guaranteed to converge at a \tcbl{linear} rate for strongly convex $f$ and sufficiently large $\mu$~\cite{dhikhojovTAC19}. \tcbl{The methods of this paper represent a second order extension of~\eqref{eq.ahu}. A comprehensive discussion of the connections between the proximal augmented Lagrangian, algorithms~\eqref{eq.MM} and~\eqref{eq.ahu}, and other commonly used techniques for solving problems~\eqref{pr} and~\eqref{pr.split} is provided in Section~\ref{sec.disc}.}

In what follows, we \tcbl{consider the continuous-time implementation of} the primal-dual \tcbl{iteration~\eqref{eq.ahu} and extend it} to incorporate second order information of $\cL_\mu(x;y)$ and thereby achieve fast convergence to high-accuracy solutions.

	\vspace*{-2ex}
\subsection{Second order updates}

	Even though Newton's method is primarily used for solving minimization problems in modern optimization, it was originally formulated as a root-finding technique and it has long been employed for finding stationary points~\cite{ortrhe00}. In~\cite{qixsun93}, a generalized Jacobian was used to extend Newton's method to semismooth problems. We \tcbl{apply} this generalization of Newton's method to $\nabla \cL_\mu(x;y)$ in order to compute the saddle point of the proximal augmented Lagrangian. The unique saddle point of $\cL_\mu(x;y)$ is given by the optimal primal-dual pair ($x^\star,y^\star$) and it thus provides the solution to~\eqref{pr}.
	
	\vsp
	
\subsubsection{Generalized Newton updates}

Let $H \DefinedAs \nabla^2 f(x)$. We use the $B$-generalized Jacobian of the proximal operator $\prox_{\mu g}$, $\bbP_B \DefinedAs \partial_B \, \prox_{\mu g}(Tx + \mu y)$, to define the set of $B$-generalized Hessians of the proximal augmented Lagrangian,
	\begin{subequations}
	\label{eq.hes-Lmu}
	\be
	\label{eq.hes-Lmu-b}
	\partial_B^2\cL_\mu
	\DefinedAs
	\left\{
	\tbt{\!\!\Hes + \frac{1}{\mu}\,T^T(I \,-\, P)T\!}{\!T^T(I \,-\, P)\!\!}{\!\!(I\,-\,P)T\!}{\!-\mu P\!\!},
	P \in \bbP_B
	\right\}
\ee
and the Clarke generalized Jacobian $\bbP_C \DefinedAs \partial_C \, \prox_{\mu g}(Tx + \mu y)$ to define the set of Clarke generalized Hessians of the proximal augmented Lagrangian,
\be
	\label{eq.hes-Lmu-c}
	\partial_C^2\cL_\mu
	\DefinedAs
	\left\{
	\tbt{\!\!\Hes + \frac{1}{\mu}\,T^T(I \,-\, P)T\!}{\!T^T(I \,-\, P)\!\!}{\!\!(I\,-\,P)T\!}{\!-\mu P\!\!},
	P \in \bbP_C
	\right\}.
\ee
\end{subequations}
Note that $\partial_B^2\cL_\mu(x;y) \subseteq \partial_C^2\cL_\mu(x;y)$ because $\bbP_B \subseteq \bbP_C$; \tcbl{$\partial_C^2\cL_\mu(x;y)$ is the convex hull of $\partial_B^2\cL_\mu(x;y)$.}

In the rest of the paper, we introduce the composite variable,
	\be
    w
    \; \DefinedAs \; 
    \obt{x^T}{y^T}^T	
    \ee
use $\cL_\mu(w)$ interchangeably with $\cL_\mu(x;y)$, and suppress the dependance of $H$ and $P$ on $w$ to reduce notational clutter. 
\tcbl{We also abuse notation and use $\partial_C^2\cL_\mu(x;y)$ to denote any member of that set.}
For simplicity of exposition, we assume that $\prox_{\mu g}$ is semismooth and state the results for the Clarke generalized Hessian~\eqref{eq.hes-Lmu-c}, i.e., $\partial^2\cL_\mu(w) = \partial_C^2\cL_\mu(w)$. As described in Remark~\ref{rem.binc} in Section~\ref{sec.diffinc}, analogous convergence results for non-semismooth $\prox_{\mu g}$ can be obtained for the $B$-generalized Hessian~\eqref{eq.hes-Lmu-b}, i.e., $\partial^2\cL_\mu(w) = \partial_B^2\cL_\mu(w)$.

We use the Clarke generalized Hessian~\eqref{eq.hes-Lmu-c} to obtain a second order update $\tilde w$ by linearizing the stationarity condition $\nabla \cL_\mu(w) = 0$ around the current iterate $w^k$,
	\beq
	\partial_C^2\cL_\mu(w^k) \, \tilde{w}^k
	\; = \; -\nabla \cL_\mu(w^k).
	\label{eq.xytilde}
	\eeq
Since $\prox_{\mu g}$ is firmly nonexpansive, $0 \preceq P \preceq I$. In Lemma~\ref{lem.sing} we use this fact to prove that the second order update $\tw$ is well-defined for any generalized Hessian~\eqref{eq.hes-Lmu} of the proximal augmented Lagrangian $\cL_\mu(x;y)$ as long as $f$ is strictly convex with $\nabla^2 f(x) \succ 0$ for all $x \in \bbR^n$.

	\vsp
	
\begin{lemma} \label{lem.sing}
Let $H \in \bbR^{n \times n}$ be symmetric positive definite, $H \succ 0$, let $P \in \bbR^{m \times m}$ be symmetric positive semidefinite with eigenvalues less than one, $0 \preceq P \preceq I$, let $T \in \bbR^{m \times n}$ be full row rank, and let $\mu > 0$. Then, the matrix
\[
	\tbt{H + \tfrac{1}{\mu} \, T^T(I - P)T}{T^T(I - P)}{(I - P)T}{-\mu P}
\]
is invertible and it has $n$ positive and $m$ negative eigenvalues.
\end{lemma}

	\vsp
	
\begin{IEEEproof}
The Haynsworth inertia additivity formula~\cite{horjoh90} implies that the inertia of matrix~\eqref{eq.hes-Lmu} is determined by the sum of the inertias of matrices,
\begin{subequations}
\be
	H \;+\; \tfrac{1}{\mu} \, T^T(I \,-\, P)T
	\label{lem.pda}
\ee
and 
\be
	-\mu P
	\, - \,
	(I - P) \, T
	\left(H + \tfrac{1}{\mu} \, T^T(I - P)T\right)^{-1}
	\!
	T^T (I - P).
	\label{lem.pdb}
\ee
\end{subequations}
Matrix~\eqref{lem.pda} is positive definite because $H \succ 0$ and both $P$ and $I - P$ are positive semidefinite. Matrix~\eqref{lem.pdb} is negative definite because the kernels of $P$ and $I - P$ have no nontrivial intersection and $T$ has full row rank.
\end{IEEEproof}

	\vsp
	
\subsubsection{Fast local convergence}
The use of generalized Newton updates for solving the nonlinear equation $G(x) = 0$ for nondifferentiable $G$ was studied in~\cite{qixsun93}. We apply this framework to the stationarity condition $\nabla \cL_\mu(w) = 0$ when $\prox_{\mu g}$ is (strongly) semismooth and show that second order updates~\eqref{eq.xytilde} converge (quadratically) superlinearly within a neighborhood of the optimal primal-dual pair.

	\vsp
	
\begin{proposition} \label{prop.quad}
Let $\prox_{\mu g}$ be (strongly) semismooth, and let $\tw^k$ be defined by~\eqref{eq.xytilde}. Then, there is a neighborhood of the optimal solution $w^\star$ in which the second order iterates
	$
	w^{k+1}
	=
	w^k
	+
	\tw^k
	$
converge (quadratically) \mbox{superlinearly to $w^\star$.}
\end{proposition}
	
\begin{IEEEproof}
Lemma~\ref{lem.sing} establishes that $\partial_C^2\,\cL_\mu(w)$ is nonsingular for any $P \in \bbP_C$. Since the gradient $\nabla \cL_\mu(w)$ of the proximal augmented Lagrangian is Lipschitz continuous by Theorem~\ref{thm.diff}, nonsingularity of $\partial_C^2\cL_\mu(w)$ and (strong) semismoothness of the proximal operator guarantee (quadratic) superlinear convergence of the iterates by~\cite[Theorem 3.2]{qixsun93}.
\end{IEEEproof}

	\vspace*{-1ex}
\section{A globally convergent differential inclusion} \label{sec.diffinc}

Since we apply a generalization of Newton's method to a saddle point problem and the second order updates are set valued, convergence to the optimal point is not immediate. Although we establish local convergence rates in Proposition~\ref{prop.quad} by leveraging the results of~\cite{qixsun93}, proof of the global convergence is more subtle and it is presented next. \tcbl{The analysis of a discrete-time algorithm is provided in Section~\ref{sec.pd}.}

To justify the development of a discrete-time algorithm based on the search direction resulting from~\eqref{eq.xytilde}, we first examine the corresponding differential inclusion,
\be
	\dot{ w}
	\;\in\;
	- (\partial_C^2\cL_\mu(w))^{-1} \, \nabla\cL_\mu(w)
	\label{eq.diffinc}
\ee
where $\partial_C^2 \cL_\mu$ is the Clarke generalized Hessian~\eqref{eq.hes-Lmu-c} of $\cL_\mu$. We assume existence of a solution and prove asymptotic stability of~\eqref{eq.diffinc} under Assumption~\ref{ass.func} and global exponential stability under an additional assumption that $f$ is strongly convex.

	\vsp

\begin{assumption} \label{ass.inc}
	Differential inclusion~\eqref{eq.diffinc} has a solution.
\end{assumption}

	\vspace*{-2ex}
\subsection{Asymptotic stability}

We first establish asymptotic stability of differential inclusion~\eqref{eq.diffinc}.

	\vsp

\begin{theorem}\label{thm.asym}
Let Assumptions~\ref{ass.func} and~\ref{ass.inc} hold and let $\prox_{\mu g}$ be semismooth. Then, differential inclusion~\eqref{eq.diffinc} is asymptotically stable. Moreover,
	\be
	V(w)
	\;\DefinedAs\;
	\tfrac{1}{2} \, \norm{\nabla \cL_\mu(w)}^2
	\label{eq.V}
\ee
provides a Lyapunov function and 
	\be
	\dot V (t) \;=\; - \, 2 V (t).
	\label{eq.gexp2}
\ee
\end{theorem}

	\vsp

\begin{IEEEproof}
Lyapunov function candidate~\eqref{eq.V} is a positive function of $w$ everywhere apart from the optimal primal-dual pair $w^\star$ where it is zero. It remains to show that $V$ is decreasing along the solutions $w(t)$ of~\eqref{eq.diffinc}, i.e., that $\dot{V}$ is strictly negative for all $w(t) \neq w^\star$,
\[
	\dot V (t)
	\;\DefinedAs\;
	\tfrac{\mathrm{d}}{\mathrm{d}t}V(w(t)) 
	\;=\;
	- \, 2V(w(t))\tcbl{.}
\]

For Lyapunov function candidates $\hat V(w)$ which are differentiable with respect to $w$, $\dot {\hat V} = \dot w^T \nabla \hat V$. Although~\eqref{eq.V} is {\em not\/} differentiable with respect to $w$, we show that $V(w(t))$ {\em is\/} differentiable  {\em along the solutions\/} of~\eqref{eq.diffinc}. Instead of employing the chain rule, we use the limit that defines the derivative,
\[
	\dot V (t)
	\, \DefinedAs \,
	\tfrac{\mathrm{d}}{\mathrm{d}t}V(w(t))
	\, = \,
	\lim_{s \, \to \, 0}
	\dfrac{V(w(t) \,+\, s \tw(t)) \,-\, V(w(t))}{s}
\]
to show that $\dot V$ exists and is negative along the solutions of~\eqref{eq.diffinc}. Here, $\tw \in - (\partial_C^2\cL_\mu(w))^{-1} \, \nabla\cL_\mu(w)$ is determined by the dynamics~\eqref{eq.diffinc}. We first introduce
\[
	h_s (t)
	\; \DefinedAs \;
	\dfrac{V(w(t) \,+\, s \tw(t)) \,-\, V(w(t))}{s}
\]
which gives $\dot V$ in the limit $s \to 0$. We then rewrite $h_s (t)$ as the limit point of a sequence of functions $\{h_{s,k} (t)\}$ so that
\be
	\dot V (t)
	\;=\;
	\lim_{s \, \to \, 0}\;
	h_s (t)
	\;=\;
	\lim_{s \, \to \, 0}\;
	\lim_{k \, \to \, \infty}\;
	h_{s,k} (t)
	\label{eq.limlim}
\ee
and use the Moore-Osgood theorem~\cite[Theorem 7.11]{rud64} to exchange the order of the limits and establish that $\dot V = -2 V$.

Let $C_g$ denote a subset of $\bbR^{n + m}$ over which $\prox_{\mu g}(Tx + \mu y)$ is differentiable (and therefore $V$ is differentiable with respect to $w$) and let $\{w_k\}$ be a sequence of points in $C_g$ that converges to $w$. We define the sequence of functions $\{ h_{s,k} (t) \}$,
\[
	h_{s,k} (t)
	\;\DefinedAs\;
	\dfrac{V(w_k(t) \,+\, s \tw_k (t)) \;-\; V(w_k(t))}{s}
\]
where $\tw_k (t) \in -( \partial_C^2 \cL_\mu (w(t)) )^{-1} \nabla\cL_\mu (w_k (t))$, as we establish below, converges to $\tw$. To employ the Moore-Osgood theorem, it remains to show that $h_{s,k} (t)$ converges pointwise (for any $k$) as $s \to 0$ and that $h_{s,k}(t)$ converges uniformly on some interval $s \in [0,\bar s]$ as $k \to \infty$.

Since $\{w_k\} \subset C_g$, $V(w_k)$ is differentiable for every $k \in \mathbb{Z}_{+}$ and $\nabla V(w_k) = \nabla^2 \cL_\mu(w_k) = \partial_C^2\cL_\mu(w_k)$. It thus follows that 
\be
\label{eq.Vdota}
	\lim_{s \, \to \, 0} \;  h_{s,k} (t) 
	\; = \;
	\tw_k^T (t)
	\, 
	\partial_C^2\cL_\mu(w_k (t))\nabla\cL_\mu(w_k (t))
\ee
pointwise (for any $k$).

We now show that the sequence $\{h_{s,k} (t)\}$ converges uniformly to $h_s (t)$ as $k \to \infty$ implying that we can exchange the order of the limits in~\eqref{eq.limlim}. Since $\prox_{\mu g}$ is semismooth, $\nabla \cL_\mu$ is semismooth. By~\eqref{eq.semi}, $\nabla \cL_\mu(w_k)$ can be written as,
\[
	\nabla \cL_\mu(w_k)
	\;=\;
	\nabla \cL_\mu(w)
	\;-\;
	\partial_C^2\cL_\mu(w_k)(w \,-\, w_k)
	\;+\;
	R_k
\]
for sufficiently large $k$, where $\norm{R_k} = o(\norm{w_k - w})$. Lemma~\ref{lem.sing} and~\cite[Proposition 3.1]{qixsun93} imply that $\partial_C^2\cL_\mu(w_k)$ is bounded within some neighborhood of $w$ and thus that $\tw_k$ can be written as
	{$
	\tw_k
	=
	\tw
	+
	\hat R_k
	$
where $\norm{\hat R_k} = O(\norm{w_k - w})$.} This implies convergence of $\tilde w_k$ to $\tilde w$ and, combined with local Lipschitz continuity of $V$ with respect to $w_k$, uniform convergence of $h_{s,k} (t)$ to $h_s (t)$ on $s \in (0,\bar s \,]$ where $\bar s >0$.

Therefore, the Moore-Osgood theorem on exchanging limits~\cite[Theorem 7.11]{rud64} in conjunction with~\eqref{eq.Vdota} imply
	\be
	\ba{rcl}
	\dot V (t)
	& \!\!\! = \!\!\! &
	\ds\lim_{s \, \to \, 0}
	h_s (t)
	\\[0.25cm]
	& \!\!\! = \!\!\! &
	\ds\lim_{s \, \to \, 0}
	\,
	\ds\lim_{k \, \to \, \infty}
	h_{s,k} (t)
	\, = \,
	\ds\lim_{k \, \to \, \infty}
	\,
	\ds\lim_{s\, \to \, 0}
	h_{s,k} (t)
	\\[0.25cm]
	& \!\!\! = \!\!\! &	
	\ds\lim_{k \, \to \, \infty}
	\tw_k^T (t) \, \partial_C^2 \cL_\mu(w_k (t))\nabla\cL_\mu(w_k (t))
	\\[0.25cm]
	& \!\!\! = \!\!\! &
	-
	\norm{\nabla \cL_\mu(w (t))}^2
	\, = \,
	-2 V(t)
	\ea
	\non
	\ee
which establishes~\eqref{eq.gexp2} and thereby completes the proof. 
\end{IEEEproof}

	\vspace*{-3ex}
\subsection{Global exponential stability}

To establish global asymptotic stability, we show that the Lyapunov function~\eqref{eq.V} is radially unbounded, and to prove exponential stability we bound it with quadratic functions. We first provide two lemmas that characterize the maps $\prox_{\mu g}$ and $\nabla f$ in terms of the spectral properties of matrices that describe the corresponding input-output relations \mbox{at given points.}

	\vsp

\begin{lemma}[Lemma 2 in~\cite{dhikhojovTAC19}]\label{lem.k}
Let $\prox_{\mu g}$: $\bbR^m \to \bbR^m$ be the proximal operator of a proper, closed, and convex function $g$. Then, for any $a,b \in \bbR^m$, there exists a symmetric matrix $D_{a,b}$ satisfying $0 \preceq D_{a,b} \preceq I$ such that 
\be
    \prox_{\mu g}(a) \,-\, \prox_{\mu g}(b) \; = \; D_{a,b} \, (a \, - \, b).
    \label{eq.lema}
    \non
\ee
\end{lemma}

	\begin{lemma} \label{lem.hatH}
Let $f$ be strongly convex with parameter $m_f$ and let its gradient $\nabla f$ be Lipschitz continuous with parameter $L_f$. Then, for any $a,b \in \bbR^n$ there exists a symmetric matrix $G_{a,b}$ satisfying $m_f I \preceq G_{a,b} \preceq L_f I $ such that 
	\be
	\nabla{f} (a) 
	\, - \, 
	\nabla{f} (b)
	\; = \; 
	G_{a,b} \, (a \, - \, b).
	\non
	\ee
\end{lemma}

\begin{IEEEproof}
Let $c \DefinedAs a - b$, $d \DefinedAs \nabla{f} (a) - \nabla{f} (b)$, $e \DefinedAs d - m_f c$, and 
	\begin{subequations}
	\begin{IEEEeqnarray}{rcl}
	\hat{G}_{a,b} 
	& ~\DefinedAs~ &
	\{ee^T/(e^Tc), \, e \neq 0; \; 0, \, \mbox{otherwise}\}
	\label{eq.hatG}
	\\[-0.1cm]
	G_{a,b}
	& ~\DefinedAs~ &
	\hat G_{a,b}
	\; + \;
	m_fI .
	\label{eq.G}
	\end{IEEEeqnarray}
	\end{subequations}
Clearly, by construction, $\hat{G}_{a,b} = \hat{G}_{a,b}^T \succeq 0$. It is also readily verified that $G_{a,b} \, c = d$ when $e^Tc \not= 0$. It thus remains to show that (i) $G_{a,b} \, c = d$ when $e^Tc = 0$; and (ii) $\hat{G}_{a,b} \preceq (L_f - m_f)I$.

(i) Since $f$ is $m_f$ strongly convex and $\nabla f$ is $L_f$ Lipschitz continuous, $h(x) \DefinedAs f(x) - \tfrac{m_f}{2} \norm{x}^2$ is convex and $\nabla h(x) = \nabla f(x) - m_f x$ is $L_f - m_f$ Lipschitz continuous. Furthermore, we have $e = \nabla h(a) - \nabla h(b)$, and~\cite[Proposition 5]{lesrecpac16} implies 
	\be
	e^Tc 
	\; \geq \;
	\tfrac{1}{L_f \, - \, m_f} \, \norm{e}^2,
	~~
	\mbox{for~all}
	~
	c \, \in \, \bbR^n.
	\label{eq.ec}
	\ee 
This shows that $e^Tc = 0$ only if $e \DefinedAs d - m_f c = 0$ and, thus, $d = m_f c = G_{a,b} \, c$ when $e^Tc = 0$. Therefore, there always exist a symmetric matrix $G_{a,b}$ such that $G_{a,b} \, c = d$.

(ii) When $e^Tc \neq 0$, $\hat{G}_{a,b}$ is a rank one matrix and its only nonzero eigenvalue is $\norm{e}^2/(e^Tc)$; this follows from $\hat{G}_{a,b} \, e = (\norm{e}^2/(e^Tc) ) \, e$. In this case, inequality~\eqref{eq.ec} implies $e^Tc > 0$ and~\eqref{eq.ec} is equivalent to $1/(e^Tc) \leq (L_f - m_f)/\norm{e}^2$. Thus, $\norm{e}^2/(e^Tc) \leq L_f - m_f$ and $\hat{G}_{a,b} \preceq (L_f - m_f) I$ when $e^Tc \neq 0$. Since $\hat{G}_{a,b} = 0$ when $e^Tc = 0$, $\hat{G}_{a,b} \preceq (L_f - m_f) I$ for all $a$ and $b$. Finally, $\hat G_{a,b} \succeq 0$ and~\eqref{eq.G} imply $m_fI \preceq G_{a,b} \preceq L_fI$.
\end{IEEEproof}

	\begin{remark} 	\label{rem.DG}
Although matrices $D_{a,b}$ and $G_{a,b}$ in Lemmas~\ref{lem.k} and~\ref{lem.hatH} depend on the operating point, their spectral properties, $0 \preceq D_{a,b} \preceq I$ and $m_f I \preceq G_{a,b} \preceq L_f I$, hold for all $a$ and $b$.
	These lemmas can be interpreted as a combination between a generalization of the mean value theorem~\cite[Theorem 5.9]{rud64} to vector-valued functions and spectral bounds on the operators $\prox_{\mu g}$: $\bbR^m \to \bbR^m$ and $\nabla f$: $\bbR^n \to \bbR^n$ arising from firm nonexpansiveness of $\prox_{\mu g}$, strong convexity of $f$, and Lipschitz continuity of $\nabla f$.
	\end{remark}
	
	\vsp
	
We now combine Lemmas~\ref{lem.k} and~\ref{lem.hatH} to establish quadratic upper and lower bounds for Lyapunov function~\eqref{eq.V} and thereby prove global exponential stability of differential inclusion~\eqref{eq.diffinc} for strongly convex $f$.
	
	\vsp
		
\begin{theorem} \label{thm.inc}
Let Assumptions~\ref{ass.func} and~\ref{ass.inc} hold, let $\prox_{\mu g}$ be semismooth, and let $f$ be strongly convex \tcbl{with parameter $m_f > 0$}. Then, differential inclusion~\eqref{eq.diffinc}  is globally exponentially stable, i.e., there exists $\kappa > 0$ such that $\norm{w(t) - w^\star} \leq \kappa \, \mre^{-t} \, \norm{w(0) - w^\star}$.
\end{theorem}

	\vsp
	
\begin{IEEEproof}
Given the assumptions, Theorem~\ref{thm.asym} establishes asymptotic stability of~\eqref{eq.diffinc} with the dissipation rate
$
	\dot V(w) \;=\; - \, 2 V(w).
$
It remains to show the existence of positive constants $\kappa_1$ and $\kappa_2$ such that Lyapunov function~\eqref{eq.V} satisfies 
\begin{IEEEeqnarray}{c}
	\tfrac{\kappa_1}{2} \, \norm{\tw}^2
		\;\leq\; 
		 V(w)
		  \;\leq\; \tfrac{\kappa_2}{2} \, \norm{\tw}^2
		  \label{eq.gexp1}
\end{IEEEeqnarray}
where $\tw \DefinedAs w - w^\star$ and $w^\star \DefinedAs (x^\star,y^\star)$ is the optimal primal-dual pair. The upper bound in~\eqref{eq.gexp1} follows from Lipschitz continuity of $\nabla \cL_\mu(w)$ (see Theorem~\ref{thm.diff}), with $\kappa_2$ determined by the Lipschitz constant of $\nabla \cL_\mu(w)$. 

To show the lower bound in~\eqref{eq.gexp1}, and thus establish radial unboundedness of $V(w)$, we construct matrices that relate $V(w)$ to $\tw$. Lemmas~\ref{lem.k} and~\ref{lem.hatH} imply the existence of symmetric matrices $D_\tw$ and $G_\tw$ such that $0 \preceq D_\tw \preceq I$,  $m_f I \preceq G_\tw \preceq L_f I$, and 
	\be
	\ba{rcl}
	\prox_{\mu g}(Tx + \mu y) - \prox_{\mu g}(Tx^\star + \mu y^\star)
	& \!\!\! = \!\!\! &
	D_\tw
	\,
	(T \tx + \mu \ty)
	\\[0.15cm]
	f(x) - f(x^\star)
	& \!\!\! = \!\!\! &
	G_\tw \, \tx.
	\ea
	\non
	\ee
As noted in Remark~\ref{rem.DG}, although $D_\tw$ and $G_\tw$ depend on the operating point, their spectral properties hold for all $\tw$.

Since $\nabla \cL_\mu(w^\star) = 0$, we can write
	$
	\nabla \cL_\mu(w)
	= 
	\nabla \cL_\mu(w)
	-
	\nabla \cL_\mu(w^\star)
	= 
	Q_\tw \tw
	$
and express Lyapunov function~\eqref{eq.V} as	
	\be
	V(w) 
	\; = \;
	\tfrac{1}{2} \,
	\tw^T Q_\tw^T \, Q_\tw \, \tw
	\non
	\ee
where
	\be
	Q_\tw 
	\; \DefinedAs \;
	\tbt{G_\tw + \tfrac{1}{\mu} \, T^T(I - D_\tw)T}{T^T(I - D_\tw)}{(I - D_\tw)T}{-\mu D_\tw}
	\non
	\ee
for some $(D_\tw,G_\tw) \in \Omega_\tw$,
\[
	\Omega_\tw
	\, \DefinedAs \,
	\left\{
	(D_\tw,G_\tw)
	~|~
	0 \, \preceq \, D_\tw \, \preceq \, I,
	~
	m_f I \, \preceq \, G_\tw \, \preceq \, L_f I
	\right\}.
\]
The set $\Omega_\tw$ is closed and bounded and the minimum eigenvalue of  $Q_\tw^T \, Q_\tw$ is a continuous function of $G_\tw$ and $D_\tw$. Thus, the extreme value theorem~\cite[Theorem 4.14]{rud64} implies that its infimum over $\Omega_\tw$,
	$
	\kappa_1
	=\;
	\inf\limits_{(D_\tw,G_\tw) \,\in\, \Omega_\tw}
	\lambda_{\min}
	\left(
	Q_\tw^T \, Q_\tw
	\right),
	$
is achieved. By Lemma~\ref{lem.sing}, $Q_\tw$ is a full rank matrix, which implies that $Q_\tw^T \, Q_\tw \succ 0$ for all $\tw$ and therefore that $\kappa_1$ is positive. Thus, $V(w) \geq \tfrac{\kappa_1}{2} \, \norm{\tw}^2$, establishing condition~\eqref{eq.gexp1}.

Condition~\eqref{eq.gexp2} and~\cite[Lemma 3.4]{kha02} imply $V(w(t)) = \mre^{-2t} V(w(0))$. It then follows from~\eqref{eq.gexp1} that
\[
	\norm{w(t) \, - \, w^\star}^2 
	\; \leq \;
	(\kappa_2/\kappa_1) \, \mre^{-2t} \, \norm{w(0) \, - \, w^\star}^2.
\]
Taking the square root completes the proof and provides an upper bound for the constant $\kappa$, $\kappa \leq \sqrt{\kappa_2/\kappa_1}$.
\end{IEEEproof}

	\vsp

\begin{remark}
The {\em rate\/} of exponential convergence established by Theorem~\ref{thm.inc} is independent of $m_f$, $L_f$, and $\mu$. This is a consequence of insensitivity of Newton-like methods to poor conditioning. In contrast, the first order primal-dual method considered in~\cite{dhikhojovTAC19} requires a sufficiently large $\mu$ for exponential convergence. In our second order primal-dual method, problem conditioning and parameter selection affect the multiplicative constant $\kappa$ but not the rate of convergence.
\end{remark}

	\vsp

	\begin{remark} \label{rem.binc}
When differential inclusion~\eqref{eq.diffinc} is defined with the $B$-generalized Hessian~\eqref{eq.hes-Lmu-b}, Theorems~\ref{thm.asym} and~\ref{thm.inc} hold even for proximal operators which are not semismooth. This follows from defining $\tw_k = \tw \in - (\partial^2_B\cL_\mu(w))^{-1}\nabla \cL_\mu(w)$ and choosing $\{w_k\} \subset C_g$ such that $\partial^2_B\cL_\mu(w) = \lim_{k \to \infty} \nabla^2\cL_\mu(w_k)$ in the proof of Theorem~\ref{thm.asym}. Such a choice of $\{w_k\}$ is possible by the definition of the $B$-subdifferential~\eqref{eq.Bsub}; since the Clarke subgradient is the convex hull of the $B$-subdifferential, it contains points outside of $\partial_B^2\cL_\mu(w)$ and thus such a sequence $\{w_k\} \subset C_g$ is not guaranteed to exist for any $\partial_C^2\cL_\mu(w)$. When defined in this manner, $\tw_k$ is constant with respect to $w_k$ and thus uniform convergence of $h_k(t,s)$ to $h(t,s)$ is immediate. 
\end{remark}

	\vspace*{-2ex}
\section{A second order primal-dual algorithm} \label{sec.pd}

An algorithm based on the second order updates~\eqref{eq.xytilde} requires step size selection to ensure global convergence. This is challenging for saddle point problems because standard notions, such as sufficient descent, cannot be applied to assess the progress of the iterates. Instead, it is necessary to identify a merit function whose minimum lies at the stationary point and whose sufficient descent can be used to evaluate progress towards the saddle point.

An approach based on discretization of differential inclusion~\eqref{eq.diffinc} and Lyapunov function~\eqref{eq.V} as a merit function leads to Algorithm~\ref{alg.v} in Appendix~\ref{app.alg}. However, such a merit function is nonconvex and nondifferentiable in general which makes the utility of backtracking (e.g., the Armijo rule) unclear. Moreover,  Algorithm~\ref{alg.v} employs a fixed penalty parameter $\mu$. {\em A priori\/} selection of this parameter is difficult and it has a large effect on \mbox{the convergence speed.}

Instead, we employ the primal-dual augmented Lagrangian introduced in~\cite{gilrobCOA12} as a merit function and incorporate an adaptive $\mu$ update. This merit function is convex in both $x$ and $y$ and it facilitates an implementation with outstanding practical performance.  Drawing upon recent advancements for constrained optimization of twice differentiable functions~\cite{armbenomhpat14,armomh15}, we show that our algorithm converges to the solution of~\eqref{pr.split}. Finally, our algorithm exhibits local (quadratic) superlinear convergence for (strongly) semismooth $\prox_{\mu g}$. 

	\vspace*{-2ex}
\subsection{Merit function}

The primal-dual augmented Lagrangian,
\be
	\non
	\cM_\mu(x,z;y,\lambda)
	\, \DefinedAs \, 
	\cL_\mu(x,z;\lambda) 
	\, + \, 
	\tfrac{1}{2\mu} \, \norm{Tx - z + \mu \, (\lambda - y)}^2
\ee
was introduced in~\cite{gilrobCOA12}, where $\lambda$ is an estimate of the optimal Lagrange multiplier $y^\star$. Following~\cite[Theorem 3.1]{gilrobCOA12}, it can be shown that the optimal primal-dual pair $(x^\star,z^\star; y^\star)$ of optimization problem~\eqref{pr.split} is a stationary point of $\cM_\mu(x,z;y,y^\star)$. Furthermore, for any fixed $\lambda$, $\cM_\mu$ is a convex function of $x$, $z$, and $y$ and it has a unique global minimizer.

In contrast to~\cite{gilrobCOA12}, we study problems in which a component of the objective function is not differentiable. As in Theorem~\ref{thm.diff}, the Moreau envelope associated with the nondifferentiable component $g$ allows us to eliminate the dependence of the primal-dual augmented Lagrangian $\cM_\mu$ on $z$,
\[
	\hat z_\mu^\star(x;y,\lambda)
	\; = \;
	\argmin\limits_z
	\,
	\cM_\mu (x, z; y, \lambda)
	\; = \;
	\prox_{\frac{\mu}{2}g}
	(Tx \,+\, \tfrac{\mu}{2}(2\lambda \,-\, y))
\]
and to express $\cM_\mu$ as a continuously differentiable function,
	\be
	\cM_\mu(x;y,\lambda)
	\; \DefinedAs \; 
	\cM_\mu(x,\hat z_\mu^\star(x;y,\lambda);y,\lambda)
	\; = \;
	\fx(x)
	\, + \,
	\mor_{\frac{\mu}{2} \gz}
	\!
	\left(Tx \,+\, \tfrac{\mu}{2}(2\lambda \,-\, y)\right)
	\, + \,
	 \tfrac{\mu}{4} \, \norm{y}^2
	 \, - \,
	 \tfrac{\mu}{2} \, \norm{\lambda}^2.
	\non
	\ee
For notational compactness, we suppress the dependence on $\lambda$ and write $\cM_\mu(w)$ when $\lambda$ is fixed. 

	\vsp

	\begin{remark}
The primal-dual augmented Lagrangian is not a Lyapunov function unless $\lambda = y^\star$. We establish convergence by minimizing $\cM_\mu(x;y,\lambda)$ over $(x;y)$ -- a convex problem -- while adaptively updating the Lagrange multiplier estimate $\lambda$.
	\end{remark}
	
	\vsp

In~\cite{gilrobCOA12,armomh15}, the authors obtain a search direction using the Hessian of the merit function,
$
    \nabla^2 \cM_\mu.
$
Instead of implementing an analogous update using generalized Hessian $\partial^2\cM_\mu$, we take advantage of the efficient inversion of $\partial^2 \cL_\mu$ (see Section~\ref{sec.eff}) to define the update
\be
	\partial^2\cL_\mu(w^k)
	\,
	\tw
	\;=\;
	-\,
	\blkdiag(I,-I) 
	\,
	\nabla \cM_{2\mu}(w^k)
	\label{eq.algxytilde}
\ee
where the identity matrices are sized conformably with the dimensions of $x$ and $y$, and
	\be
	\nabla\cM_{2\mu}(w)
	= 
	\left[
	\ba{c}
	\gr(x) 
	+ 
	T^T 
	\nabla
	\mor_{\mu \gz}
	\!
	\left(Tx + \mu (2\lambda - y)\right)
	\\
	\mu (y - 
	\nabla
	\mor_{\mu \gz}
	\!
	\left(Tx + \mu (2\lambda - y)\right))
	\ea
	\right].
	\label{eq.gradMer}
	\ee
Multiplication by $\blkdiag(I,-I)$ is used to ensure descent in the dual direction and $\cM_{2\mu}$ is employed because $\hat z^\star_\mu (x;y,\lambda)$ is determined by the proximal operator associated with $(\mu/2) g$. When $\lambda = y$, $\nabla_x \cM_{2\mu} = \nabla_x \cL_{\mu}$, $\nabla_y \cM_{2\mu} = - \nabla_y \cL_{\mu}$, and~\eqref{eq.algxytilde} becomes equivalent to the second order update~\eqref{eq.xytilde}.

	\vsp

\begin{lemma}
\label{lem.descent}
Let $\tw$ solve~\eqref{eq.algxytilde}. Then, for the fixed value of the Lagrange multiplier estimate $\lambda$ and any $\sigma \in (0,1]$, 
\be
	d
	\; \DefinedAs \;
	(1 \, - \, \sigma) \, \tw
	\;-\;
	 \sigma \, \nabla \cM_{2\mu}(w,\lambda)
\label{eq.d}
\ee
is a descent direction of the merit function $\cM_{2\mu}(w,\lambda)$.
\end{lemma}

	\vsp

\begin{IEEEproof}
By multiplying~\eqref{eq.algxytilde} with the nonsingular matrix
\be
	\Pi
	\;\DefinedAs\;
	\tbt{I}{-\,\tfrac{1}{\mu}\,T^T}{0}{I}
	\label{eq.Pi}
\ee
we can express it as
\be
	\label{eq.search}
	\tbt{\!\!\Hes\!\!}{\!\!T^T\!\!}{\!\!(I \,-\, P)T\!\!}{\!\!-\mu P\!\!}
	\tbo{\tx}{\ty}
	= 
	\tbo{- (\gr(x) + T^Ty)}{\nabla_y \cM_{2 \mu} (x; y, \lambda)}
\ee
where $H \DefinedAs \nabla^2 f (x) \succ 0$. Using~\eqref{eq.gradMer} and~\eqref{eq.search}, $\nabla \cM_{2\mu}(w)$ can be expressed as,
\[
	\nabla \cM_{2\mu}(w)
	= 
	\tbo{-(\Hes + \tfrac{1}{\mu}T^T(I-P)T ) \tx \, - \, T^T(I - P) \ty}
	{ (I\,-\,P)T\tx \;-\; \mu P \ty}.
	\non
\]
Thus, 
\[
	{\tw}^T {\nabla \cM_{2\mu}(w)}
	\, = \,
	 - \tx^T (\Hes + \tfrac{1}{\mu}T^T(I\,-\,P)T) \tx  \, - \, \mu\ty^TP\ty
\]
is negative semidefinite, and the inner product
\[
    \ba{c}
	{d}^T {\nabla \cM_{2\mu}(w)}
	\, = \,
      (1-\sigma) 
      \, 
      {\tw}^T {\nabla \cM_{2\mu}(w)}
        \,-\, 
        \sigma\norm{\nabla\cM_{2\mu} (w)}^2
    \ea
\]
is negative definite when $\nabla \cM_{2\mu}$ is nonzero.
\end{IEEEproof}

	\vspace*{-2ex}
\subsection{Second order primal-dual algorithm}
\label{sec.alg}

We now develop a customized algorithm that alternates between minimizing the merit function $\cM_\mu(x;y,\lambda)$ over $(x;y)$ and updating $\lambda$. Near the optimal solution, the algorithm approaches second order updates~\eqref{eq.xytilde} with unit step size, leading to local (quadratic) superlinear convergence for (strongly) semismooth $\prox_{\mu g}$.

Our approach builds on the sequential quadratic programming method described in~\cite{gilrobCOA12,armbenomhpat14,armomh15} and it uses the primal-dual augmented Lagrangian as a merit function to assess progress of iterates to the optimal solution. Inspired by~\cite{delfackan96}, we ensure sufficient progress with damped second order updates.

The following two quantities
	\[
	\ba{rcl}
	r
	& \!\!\! \DefinedAs \!\!\! &
	Tx \, - \, \prox_{\mu \gz}
	\!
	\left(Tx \, + \, \mu y \right)
	\\[0.1cm]
	s
	& \!\!\! \DefinedAs \!\!\! &
	Tx \, - \, \prox_{\mu \gz}
	\!
	\left(Tx \, + \, \mu \, (2\lambda \, - \, y) \right)
	\ea
	\]
appear in the proof of global convergence. Note that $r$ is the primal residual of optimization problem~\eqref{pr.split} and that $\nabla\cM_{2\mu}$ can be equivalently expressed as
	\be
	\nabla\cM_{2\mu}(w)
	\;=\;
	\left[
	\ba{c}
	\gr(x) \, + \, \tfrac{1}{\mu} \, T^T (s + \mu(2\lambda - y))
	\\
	-(s \, + \, 2\mu(\lambda - y))
	\ea
	\right].
	\label{eq.gradMs}
	\ee

	\vsp
	
\subsubsection{Global convergence}

We now establish global convergence of Algorithm~\ref{alg.sqp} under \tcbl{either an assumption on $g$ being finite valued or} an assumption that the sequence of gradients generated by the algorithm, $\nabla f(x^k)$, is bounded. The \tcbl{latter} assumption is standard for augmented Lagrangian based methods~\cite{armomh15,congoutoi91} and it does not lead to a loss of generality when $f$ is strongly convex. \tcbl{We first show a small Lemma to facilitate the proof.}

	\vsp
\begin{lemma} \label{lem.proxgfinite}
	\tcbl{Let $g$ be finite everywhere. Then, as $\mu \to 0$, $\prox_{\mu g}(v) \to v$.}
\end{lemma}

	\vsp

\begin{IEEEproof}
	\tcbl{The optimality condition of~\eqref{eq.prox} yields the relation $\prox_{\mu g}(v) = v - \mu\,\partial g (\prox_{\mu g}(v))$. Since $g$ is finite everywhere, $\partial g$ is finite everywhere. Thus, as $\mu \to 0$, $\mu\,\partial g (\prox_{\mu g}(v)) \to 0$ and $\prox_{\mu g}(v) \to v$. }
\end{IEEEproof}

	\vsp
	
\begin{theorem} \label{thm.conv}
Let Assumption~\ref{ass.func} hold and let \tcbl{either $g$ be finite everywhere or} the sequence $\{\nabla f(x^k)\}$ resulting from Algorithm~\ref{alg.sqp} be bounded. Then, the sequence of iterates $\left\{w^k\right\}$ converges to the optimal primal-dual point of problem~\eqref{pr.split} and the Lagrange multiplier estimates $\{ \lambda^k \}$ converge to the optimal Lagrange multiplier.
\end{theorem}

	\vsp

\begin{IEEEproof}
Since ${\cal V}_{2\mu}(w,\lambda)$ is convex in $w$ for any fixed $\lambda$, condition~\eqref{alg.req2} in Algorithm~\ref{alg.sqp} will be satisfied after finite number of iterations. Combining~\eqref{alg.req2} and~\eqref{eq.gradMs} shows that $s^k + 2\mu^k(\lambda^k - y^k) \to 0$ and $\gr(x^k) + \tfrac{1}{\mu^k} \, T^T (s^k + \mu(2\lambda^k - y^k))  \to 0$. Together, these statements imply that the dual residual $\nabla f(x^k) + T^Ty^k$ of~\eqref{pr.split} converges to zero.

To show that the primal residual $r^k$ converges to zero, we first show that $s^k \to 0$. If Step $2a$ in Algorithm~\ref{alg.sqp} is executed infinitely often, $s^k \to 0$ since it satisfies~\eqref{alg.req1} at every iteration and $\eta \in (0,1)$. If Step $2a$ is executed finitely often, there is $k_0$ after which $\lambda^k = \lambda^{k_0}$. 

\tcbl{The second row of~\eqref{alg.req2} and~\eqref{eq.gradMs} imply that $s^k \to 2\mu^k y^k$ since $\mu^k \to 0$ and $\lambda^k = \lambda^{k_0}$. However, if $g$ is finite everywhere, Lemma~\ref{lem.proxgfinite} and the definition of $s$ imply that $s^k \to \mu^k y^k$. Therefore, it must be that $s^k \to 0$.}

\tcbl{If $g$ is not finite everywhere, $s^k \to 0$ is implied by boundedness of  $\nabla f(x^k)$. In particular,} by adding and subtracting $2\mu^k \nabla f(x^k) + T^Ts^k + 4\mu^kT^T(\lambda^{k_0} - y^k)$ and rearranging terms, we can write
	\[
    T^Ts^k
    \; = \;
    2\mu^k(\nabla f(x^k) \,+\, \tfrac{1}{\mu^k}T^T(s^k + \mu^k(2\lambda^{k_0} - y^k)))
    \, - \,
    2\mu^k \nabla f(x^k)
     \,-\,
    T^T(s^k + 2\mu^k(\lambda^{k_0} - y^k))
    \, - \,
    2\mu^k T^T\lambda^{k_0}.
	\]	
\noindent
Taking the norm of each side and applying {the triangle inequality,}~\eqref{alg.req2} and~\eqref{eq.gradMs} yields
\be
    \norm{T^Ts^k}
    \, \leq \,
    {2
    \mu^k}
    \epsilon^k
     + 
    2\mu^k 
    \norm{\nabla f(x^k)}
    +
    \norm{T^T}\epsilon_k
     + 
   2\mu^k
    \norm{T^T\lambda^{k_0}}.
\label{thm.ineq}
\ee
This inequality implies that $T^Ts^k \to 0$ because $\nabla f(x^k)$ is bounded, $\epsilon^k \to 0$, and $\mu^k \to 0$. Since $T$ has full row rank, $T^T$ has full column rank and it follows that $s^k \to 0$.

Substituting $s^k \to 0$ and $\nabla f(x^k) + T^Ty^k \to 0$ into the first row of~\eqref{eq.gradMs} and applying~\eqref{alg.req2} implies $\lambda^k \to y^k$. Thus, $s^k \to r^k$, implying that the iterates asymptotically drive the primal residual $r^k$ to zero, thereby completing the proof.
\end{IEEEproof}

	\vsp

	\begin{remark}
Despite the assumption that $\{\nabla f(x^k)\}$ is bounded, Theorem~\ref{thm.conv} can be used to ensure global convergence whenever $f$ is strongly convex. We show in Lemma~\ref{lem.bdset} in Appendix~\ref{sec.appf} that a bounded set $\cC_f$ containing the optimal point can be identified {\em a priori\/}. One can thus artificially bound~$\nabla f(x)$ for all $x \not\in \cC_f$ to satisfy the conditions of Theorem~\ref{thm.conv} and guarantee global convergence \mbox{to the solution of~\eqref{pr}.}
	\end{remark}
	
	\vsp
	\begin{algorithm}
\caption{Second order primal-dual algorithm for nonsmooth composite optimization.}
\label{alg.sqp}
\begin{algorithmic}
\STATE \textbf{input:} Initial point $w^0 = (x^0,y^0)$, \tcbl{tolerance} ${\eta} \in (0,1)$, \tcbl{initial stepsize $\alpha \in (0,1)$,} Armijo parameter $\beta \in (0,1)$, \tcbl{$\mu$ update parameters} $\tau_a,\tau_b \in (0,1)$, and \tcbl{decreasing intermediate tolerances} $\{\epsilon\}^k \geq 0$ such that $\epsilon^k \rightarrow 0$.
\STATE \textbf{initialize:} Set $\lambda^0 = y^0$.
\vspace*{0.15cm}
\\
~\,\quad \textbf{Step $1$:} If
\be
    \norm{s^k}
    \;\leq\;
    {\eta}
    \norm{s^{k-1}}
\label{alg.req1}
\ee
~\,\quad \quad\quad go to Step $2a$. If not, go to Step $2b$.

~\,\quad \textbf{Step $2a$:} Set
\[
    \ba{rclrcl}
    \mu^{k+1}
    &\!\!\! = \!\!\!&
    \tau_a \mu^k,
    &
    \hspace{0.5cm}
    \lambda^{k+1}
    &\!\!\! = \!\!\!&
    y^k
    \ea
\]
~\,\quad \textbf{Step $2b$:} Set
\[
    \ba{rclrcl}
    \mu^{k+1}
    &\!\!\! = \!\!\!&
    \tau_b \mu^k,
    &
    \hspace{0.5cm}
    \lambda^{k+1}
    &\!\!\! = \!\!\!&
    \lambda^k
    \ea
\]	
~\,\quad \textbf{Step $3$:} \tcbl{Perform a sequence of inner iterations $w^{k,j}$ until some iterate $J$ such that $w^{k+1} = w^{k,J}$ satisfies}
\be
\label{alg.req2}
    \norm{\nabla \cM_{2\mu^{k+1}}(w^{k+1},\lambda^{k+1})}
    \;\leq\;
    \epsilon^k.
\ee
\tcbl{Setting $w^{k,1} = w^k$, ensuing inner iterations are of the form
\[
    w^{k,j+1} 
    \;=\;
    w^{k,j}
    \,+\,
    \alpha^l d^{k,j},
\]
where $l$ is the smallest integer that satisfies the Armijo rule~\cite{nocwri06} and the search direction $d^{k,j}$ is obtained from~\eqref{eq.d}--\eqref{eq.search} with $(w^{k,j},\lambda^k)$ and choosing $\sigma$ in~\eqref{eq.d} as}
\begin{subequations}
\begin{IEEEeqnarray}{lr}
	\sigma \, = \, 0,
	&
	~~
	\dfrac{(\tw^{k})^T \nabla\cM_{2\mu^{k+1}}(w^{k})}{\norm{\nabla \cM_{2\mu^{k+1}}(w^{k})}^2} \;\leq\; -\beta,  \label{eq.thmsig}
	\\
	\sigma \, \in \, (0,1],
	&
	\mbox{otherwise.}
\end{IEEEeqnarray}
\end{subequations}
\end{algorithmic}
\end{algorithm}

\subsubsection{Asymptotic convergence rate}

The invertibility of the generalized Hessian $\partial^2 \cL_\mu(w)$ allows us to establish local convergence rates for the second order updates~\eqref{eq.xytilde} when $\prox_{\mu g}$ is (strongly) semismooth.
	
We now show that the updates in Algorithm~\ref{alg.sqp} are equivalent to the second order updates~\eqref{eq.xytilde} as $k \to \infty$. Thus, if $\prox_{\mu g}$ is (strongly) semismooth, the sequence of iterates generated by Algorithm~\ref{alg.sqp} converges (quadratically) superlinearly to the optimal point in some neighborhood of it.

	\vsp

\begin{theorem}
Let the conditions of Theorem~\ref{thm.conv} hold, let $\prox_{\mu g}$ be (strongly) semismooth, and let $\epsilon^k$ be such that $\norm{w^k - w^\star} = O(\epsilon^k)$. Then, in a neighborhood of the optimal point $w^\star$, $w^k$ converges (quadratically) superlinearly to $w^\star$.
\end{theorem}

	\vsp
	
\begin{IEEEproof}
	From Theorem~\ref{thm.conv}, $\lambda^k \to y^k$ and thus $\nabla \cM_{2\mu}(w^k) \to \nabla \cL_\mu(w^k)$. Descent of the Lyapunov function in Theorem~\ref{thm.asym} therefore implies that the update in Step $3$ of Algorithm~\ref{alg.sqp} is given by~\eqref{eq.thmsig}, which is equivalent to~\eqref{eq.xytilde} because $\lambda = y$. The assumption on $\{\epsilon^k\}$ in conjunction with Proposition~\ref{prop.quad}, and~\cite[Theorem 3.2]{fac95} imply that this update asymptotically satisfies~\eqref{alg.req2} in one iteration with a unit step size. Therefore, Step $3$ reduces to~\eqref{eq.xytilde} in some neighborhood of the optimal solution and Proposition~\ref{prop.quad} implies that $w^k$ converges to $w^\star$ (quadratically) superlinearly.
\end{IEEEproof}

	\vsp

\subsubsection{Efficient computation of the Newton direction}
\label{sec.eff}

When $g$ is (block) separable, the matrix $P$ in~\eqref{eq.hes-Lmu} is (block) diagonal. We next demonstrate that the solution to~\eqref{eq.search} can be efficiently computed when $T = I$ and $P$ is a sparse diagonal matrix whose entries are either $0$ or $1$. The extensions to a low rank $P$, to a $P$ with entries between $0$ and $1$, or to a general diagonal $T$ follow from similar arguments.

These conditions occur, for example, when $\gz(z) = \gamma \| z \|_1$. The matrix $P$ is sparse when $\prox_{\mu g}(x + \mu y) = \cS_{\gamma\mu}(x + \mu y)$ is sparse. Larger values of $\gamma$ are more likely to produce a sequence of iterates $w^k$ for which $P$ is sparse and thus the second order search directions~\eqref{eq.search} are cheaper to compute.

We can write~\eqref{eq.search} as
\begin{subequations}
\be
    \tbt{H}{I}{I-P}{-\mu P}
    \tbo{\tx}{\ty}
    \;=\;
    \tbo{\vartheta}{\theta}
    \label{eq.stilde}
\ee
permute it according to the entries of $P$ which are $1$ and $0$, respectively, and partition the matrices $\Hes$, $P$, and $I - P$ conformably such that
\[
	\ba{l}
	\Hes
	\;=\;
	\tbt{\Hes_{11}}{\Hes_{12}}{\Hes_{12}^T}{\Hes_{22}},
	~
	P
	\;=\;
	\tbt{I}{}{}{0}.
	\ea
\]
Let $v$ denote either the primal variable $x$ or the dual variable $y$. We use $v_1$ to denote the subvector of $v$ corresponding to the entries of $P$ which are equal to $1$ and $v_2$ to denote the subvector corresponding to the zero diagonal entries of $P$. 
	
Note that $(I - P)v = 0$ when $v_2 = 0$ and $Pv = 0$ when $v_1 = 0$. As a result, $\tx_2$ and $\ty_1$ are explicitly determined by the bottom row of the system of equations~\eqref{eq.stilde},
\be
	\label{eq.xbpyp}
	\tbt{0}{-\mu I}{I}{0}
	\tbo{\tx_2}{\ty_1}
	\;=\;
	\tbo{\theta_1}{\theta_2}
\ee
Substitution of the subvectors $\tx_2$ and $\ty_1$ into~\eqref{eq.stilde} yields,
\be
	\label{eq.xp}
	H_{11}\tx_1
	\;=\;
	\vartheta_1 \,+\, H_{12}\tx_2 \,+\, \ty_1
\ee
which must be solved for $\tx_1$. Finally, the computation of $\ty_2$ requires only matrix-vector products,
\be
	\label{eq.ybp}
	\ty_2
	\;=\;
	- \left( \vartheta_2 \,+\, \Hes_{21}\tx_1 \,+\, \Hes_{22}\tx_2\right).
\ee
\end{subequations}
Thus, the major computational burden in solving~\eqref{eq.search} lies in performing a Cholesky factorization to solve~\eqref{eq.xp}, where $\Hes_{11}$ is a matrix of a much smaller size than $\Hes$.

	\vspace*{-1ex}
\section{Computational experiments}
\label{sec.ex}

In this section, we illustrate the merits and the effectiveness of our approach. \tcbl{We use an $\ell_1$-regularized model predictive control (MPC) problem with $g = \gamma \|x\|_1$ to demonstrate a speedup for larger values of the regularization parameter $\gamma$. We then consider a controller synthesis problem for a spatially-invariant system with $T \neq I$. In both examples, we choose $\epsilon^k = \norm{\nabla\cV_{2\mu^0}(w^0,\lambda^0)}/k$ and set $(\eta,\alpha,\beta,\tau_a,\tau_b,\mu^0,\sigma)$ to $(0.8,0.5,0.001,0.6,0.6,100,0.001)$.}

	\vspace*{-2ex}
\subsection{$\ell_1$-regularized model predictive control problem}

\begin{figure}
\begin{center}
	\[
	\ba{cc}
	{\footnotesize \gamma
	\;=\;
	0.15\gamma_{\max}}
    \!\!\!\!\!
    &
    \!\!\!\!\!
	{\footnotesize \gamma
	\;=\;
	0.85\gamma_{\max}}
	\\
	  \centering
  \begin{tabular}{rc}
      \!\!\!\!\!
    \rotatebox{90}{\footnotesize \quad \quad \quad \quad  $\norm{x \,-\, x^\star}$}
    \!\!\!\!\!
    &
    \!\!\!\!\!
    {	\includegraphics[width= 0.48\columnwidth]{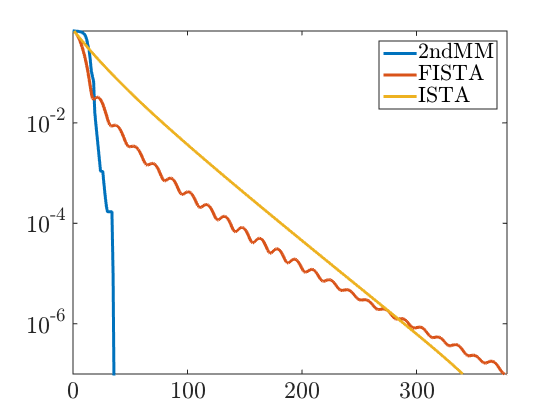}}
    \\
    \!\!\!\!\!
    &
    \!\!\!\!\!
    {\footnotesize iteration}
    \end{tabular}
    \!\!\!\!\!    \!\!\!\!\!
    &
    \!\!\!\!\!    \!\!\!\!\!
\begin{tabular}{rc}
    \!\!\!\!\!
    &
    \!\!\!\!\!
    {	\includegraphics[width= 0.48\columnwidth]{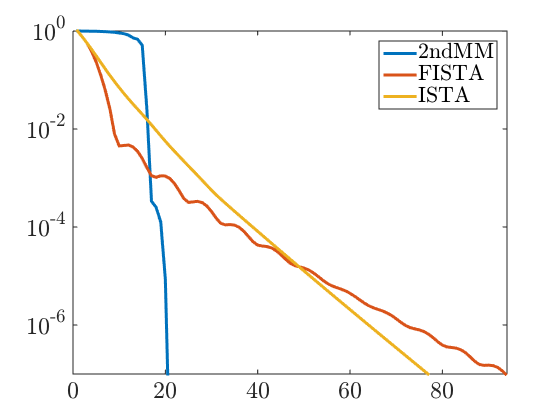}}
    \\
    \!\!\!\!\!
    &
    \!\!\!\!\!
    {\footnotesize iteration}
    \end{tabular}
	\\
	  \centering
  \begin{tabular}{rc}
      \!\!\!\!\!
    \rotatebox{90}{\footnotesize \quad \quad \quad \quad  $\norm{x \,-\, x^\star}$}
    \!\!\!\!\!
    &
    \!\!\!\!\!
    {	\includegraphics[width= 0.48\columnwidth]{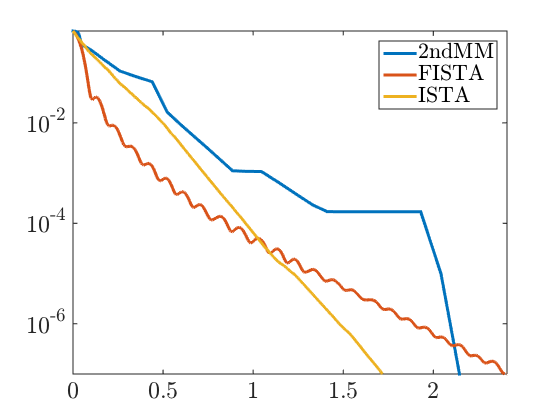}}
    \\
    \!\!\!\!\!
    &
    \!\!\!\!\!
    {\footnotesize time (s)}
    \end{tabular}
    \!\!\!\!\!    \!\!\!\!\!
    &
    \!\!\!\!\!    \!\!\!\!\!
\begin{tabular}{c}
    \!\!\!\!\!
    {	\includegraphics[width= 0.48\columnwidth]{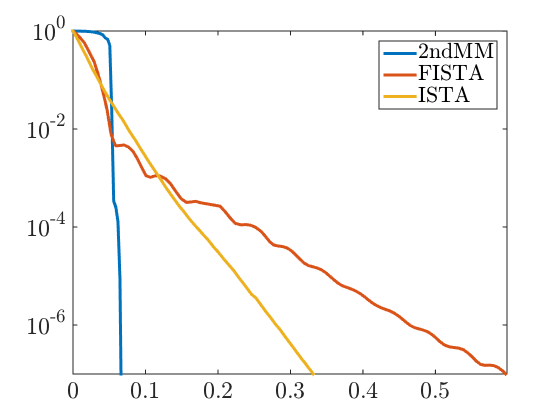}}
    \\
    \!\!\!\!\!
    {\footnotesize time (s)}
    \end{tabular}
	\ea
	\]
\end{center}
\caption{Distance from optimal solution as a function of the iteration number and solve time when solving LASSO for two values of $\gamma$ using ISTA, FISTA, and our algorithm ($2$ndMM).}
\label{fig.lassoista}
\end{figure}

\begin{figure}
\begin{center}
\[
	\begin{tabular}{rc}
    \!\!\!\!\!
        \rotatebox{90}{\footnotesize \quad Computation time (s)}
    \!\!\!\!\!
    &
    \!\!\!\!\!
    {	\includegraphics[width= 0.45\columnwidth]{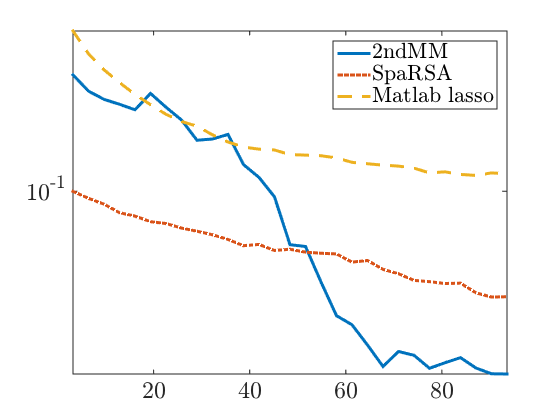}}
        \\
    \!\!\!\!\!
    \!\!\!\!\!
    &
    \!\!\!\!\!
    {\footnotesize Percent of $\gamma_{\max}$}
    \end{tabular}
\]
\end{center}
    \caption{Solve times for LASSO with $n = 1000$ obtained using ISTA, FISTA, and our algorithm ($2$ndMM) as a function of the sparsity-promoting parameter $\gamma$.}
    \label{fig.lassogam}
\end{figure}

\begin{figure}
\begin{center}
\vspace{-.5cm}
\[
  \begin{tabular}{rcrc}
    \!\!\!\!\!
    &
    \!\!\!\!\!
    {\footnotesize $\gamma \;=\; 0.15\gamma_{\max}$}
    \!\!\!\!\!
    &
    \!\!\!\!\!
    \!\!\!\!\!
    &
    \!\!\!\!\!
    {\footnotesize $\gamma \;=\; 0.85\gamma_{\max}$}
    \\
    \!\!\!\!\!
        \rotatebox{90}{\footnotesize   \quad\; Computation time (s)}
    \!\!\!\!\!
    &
        \!\!\!\!\!
    {	\includegraphics[width= 0.48\columnwidth]{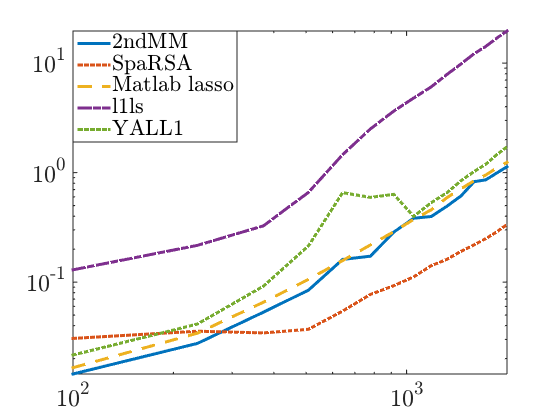}}
    \!\!\!\!\!    \!\!\!\!\!
    &
    \!\!\!\!\!    \!\!\!\!\!
    \!\!\!\!\!    \!\!\!\!\!
    &
    \!\!\!\!\!    \!\!\!\!\!
    {	\includegraphics[width= 0.48\columnwidth]{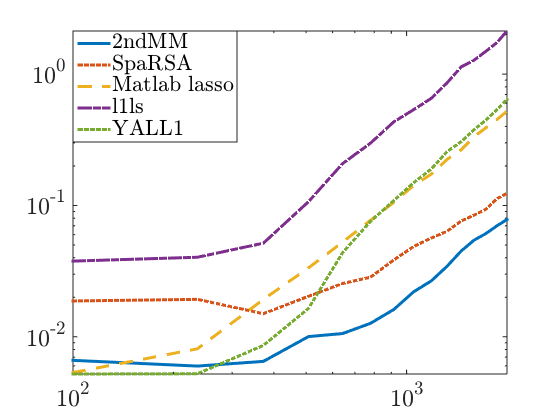}}
        \\
    \!\!\!\!\!
    &
    \!\!\!\!\!
    {\footnotesize $n$}
    \!\!\!\!\!
    &
    \!\!\!\!\!
    \!\!\!\!\!
    &
    \!\!\!\!\!
    {\footnotesize $n$}
    \end{tabular}
\]
\end{center}
\caption{Comparison of our algorithm (2ndMM) with state-of-the-art methods for LASSO with problem dimension varying from $n = 100$ to $2000$.}
\label{fig.comp}
\end{figure}

\tcbl{In the MPC framework, control inputs are determined by solving a finite horizon optimization problem. In order to minimize the number of actuators, typical quadratic penalties on state and control signals were augmented with an $\ell_1$ penalty on the control input in~\cite{galmac12}. This approach is of interest for over-actuated systems and it has been applied to naval guidance~\cite{galmac12} as well as spacecraft rendezvous and capture control~\cite{hargalmac13}. Our algorithm is well-suited to these MPC problems because they must be solved quickly and to high-accuracy to support run-time constraints and ensure the validity of the solution.}

\tcbl{For a discrete time linear system, MPC is determined by the solution to an $\ell_1$-regularized quadratic objective function over a finite time interval $t \in \{0, \dots, T_f\}$,
\be
    \begin{array}{rl}
        \minimize\limits_{\psi, \, u}
        &
        \ds{\sum\limits_{t \, = \, 0}^{T_f}}
        \left(
        \psi_t^TQ_t\psi_t
        \;+\;
        u_t^TR_tu_t
        \;+\;
        \gamma
        \,
        \| u_t \|_1
        \right)
        \\[0.35cm]
        \subject
        &
        \psi_{t+1}
        \;=\;
        A
        \,
        \psi_t
        \;+\;
        B
        \,
        u_t,
        ~~
        \psi_0
        \;=\;
        \hat \psi_0.
    \end{array}
    \label{eq.mpc}
\ee
Here, $\psi_t$ is the state, $u_t$ is the control input, $Q_t \succeq 0$ and $R_t \succ 0$ are the state and control weights, and $\hat \psi_0$ is the initial condition. Problem~\eqref{eq.mpc} can be written as an $\ell_1$-regularized quadratic problem,
\[
    \minimize\limits_x
    \left(
    (\bar \psi_0 + \bar Ax)^T \bar Q \bar A (\bar \psi_0 + \bar A x)
    \, + \,
    x^T\bar R x
    \, + \,
    \gamma \norm{x}_1
    \right)
\]
where $x \DefinedAs [u_0^T ~\cdots~ u_{T_f}^T]^T$, and
\[
    \ba{rcl}
        \bar Q
        &\!\!=\!\!&
        \thbth{Q_0}{}{}
                {}{\ddots}{}
                {}{}{Q_{T_f}},
        ~~
        \bar R
        \;=\;
        \thbth{R_0}{}{}
                {}{\ddots}{}
                {}{}{R_{T_f}}
        \\
        \bar A
        &\!\!=\!\!&
        \left[
            \ba{cccc}
            {0}&{0}&{0}&{0} \\
            {B}&{0}&{0}&{0} \\
            {\vdots}&{\ddots}&{0}&{0} \\
            {A^{T_f-1}B}&{\dots}&{B}&{0}
            \ea
        \right],
        ~~
         \bar \psi_0
        \;=\;
        \fbo{\hat{\psi}_0}{A \hat{\psi}_0}{\vdots}{A^{T_f} \hat{\psi}_0}.
    \ea
\]
Furthermore, completion of squares can be used to obtain a special form of the canonical LASSO problem,
\be
    \minimize\limits_x
    ~\,
    \tfrac{1}{2} \, \norm{\tcbl{F}x \,-\, b}^2
    \;+\;
    \gamma
    \,
    \norm{x}_1.
    \label{eq.lasso}
\ee
Since $\bar R \succ 0$, the instances of~\eqref{eq.lasso} that arise from an MPC problem~\eqref{eq.mpc} {\em are strongly convex}. This is in contrast to compressive sensing applications, where it is desired to find a sparse solution to an underdetermined system of equations.} 

\tcbl{We note that MPC formalism allows for inclusion of  additional penalty functions, time varying dynamics, as well as constraints on states and control inputs. Our algorithm can be useful for many such problems; its suitability for a particular instance depends on whether it can be cast in the form~\eqref{pr.split} where $g$ has an easily computable proximal operator.}

As described in Section~\ref{sec.back_back}, the proximal operator \tcbl{associated with the $\ell_1$ norm} is given by soft-thresholding $\cS_{\gamma \mu}$, the Moreau envelope is the Huber function, and its gradient is the saturation function. Thus, $P \in \bbP$ is diagonal and $P_{ii}$ is $0$ when $| x_i + \mu y_i | < \gamma \mu$, $1$ outside this interval, and between $0$ and $1$ on the boundary. Larger values of the regularization parameter $\gamma$ induce sparser solutions for which one can expect a sparser sequence of iterates. Note that we require strong convexity of the least squares penalty; i.e., that \tcbl{$F^TF$} is positive definite.
	
In Fig.~\ref{fig.lassoista}, we show the distance of the iterates from the optimal for the standard proximal gradient algorithm ISTA, its accelerated version FISTA, and our customized second order primal-dual algorithm for a problem where \tcbl{$F^TF$} has condition number $3.26 \times 10^4$.
\tcbl{We motivate problem~\eqref{eq.lasso} with MPC problem~\eqref{eq.mpc} because its control effort penalty results in a positive definite $F^TF$; this is contrast to many sparse inverse problems which are overdetermined, meaning that $F$ is a fat matrix and $F^TF$ is rank deficient. However, despite drawing inspiration from~\eqref{eq.mpc}, we consider a {\em generic\/} random matrix $F$ and random vector $b$ in our computational experiments in the interest of simplicity so that the computational trends we observe can be confidently ascribed to our algorithm instead of to the problem structure.} We plot distance from the optimal point as a function of both iteration number and solve time. Although our method always requires much fewer iterations, it is most effective when $\gamma$ is large. In this case the most computationally demanding step~\eqref{eq.xp} required to determine the second order search direction~\eqref{eq.search} involves a smaller matrix inversion; see Section~\ref{sec.eff} for details. In Fig.~\ref{fig.lassogam}, we show the solve times for $n=1000$ as the sparsity-promoting parameter $\gamma$ ranges from $0$ to $\gamma_{\max} = {\norm{F^T b}_\infty}$, where $\gamma_{\max}$ yields a zero solution. All numerical experiments consist of \mbox{$20$ averaged trials.}
	
In Fig.~\ref{fig.comp}, we compare the performance of our algorithm with the LASSO function in Matlab (a coordinate descent method \cite{frihastib10}), SpaRSA~\cite{wrinowfig09}, an interior point method~\cite{kimkohlusboygor07}, and YALL1~\cite{yanzha11}. Problem instances were randomly generated with $\tcbl{F} \in \bbR^{m \times n}$, $n$ ranging from $100$ to $2000$, $m = 3n$, and $\gamma = 0.15\gamma_{\max}$ or $0.85\gamma_{\max}$. The solve times and scaling of our algorithm is competitive with these state-of-the-art methods. For larger values of $\gamma$, the second order search direction~\eqref{eq.search} is cheaper to compute and our algorithm is the fastest.

	\vspace*{-2ex}
\subsection{Distributed control of a spatially-invariant system}

We now apply our algorithm to a structured control design problem aimed at balancing closed-loop $\cH_2$ performance with spatial support of a state-feedback controller. Following the problem formulation of~\cite{linfarjovTAC13admm}, ADMM was used in~\cite{wujovSCL17} to design sparse feedback gains for spatially-invariant systems. Herein, we demonstrate that our algorithm provides significant computational advantage over both the ADMM algorithm and a proximal Newton scheme.

\subsubsection{Spatially-invariant systems}
Let us consider 
\be
	\ba{rcl}
	\dot \psi
	&\!\!=\!\!&
	A\,\psi
	\;+\;
	u
	\;+\;
	d
	\\[0.1cm]
	\zeta
	&\!\!=\!\!&
	\tbo{Q^{1/2}\,\psi}
	{R^{1/2}\,u}
	\ea
	\label{eq.spinv}
\ee
where $\psi$, $u$, $d$, and $\zeta$ are the system state, control input, white stochastic disturbance, and performance output and $A$, $Q \succeq 0$, and $R \succ 0$ are $n \times n$ circulant matrices. Such systems evolve over a discrete spatially-periodic domain; they can be used to model spatially-invariant vehicular platoons~\cite{bamjovmitpat12} and can result from a spatial discretization of fluid flows~\cite{jovbamJFM05}.

Any circulant matrix can be diagonalized via the discrete Fourier transform (DFT).  Thus, the coordinate transformation $\psi \DefinedAs T \hat{\psi}$, $u \DefinedAs T \hat{u}$, $d \DefinedAs T \hat{d}$, where $T^{-1}$ is the DFT matrix, brings the state equation in~\eqref{eq.spinv} into,
\be
	\dot{\hat{\psi}}
	\; = \;
	\hat{A} \, \hat{\psi}
	\;+\;
	\hat{u}
	\;+\;
	\hat{d}.
	\label{eq.spinv-hat}
\ee
Here, $\hat{A} \DefinedAs T^{-1} A T$ is a diagonal matrix whose main diagonal $\hat{a}$ is determined by the DFT of the first row $a$ \mbox{of the matrix $A$,} 
	\be
	\hat{a}_k
	\; \DefinedAs \;
	\sum_{i \, = \, 0}^{n-1}
	\, 
	a_i 
	\,
	\mre^{- \mrj \tfrac{2 \pi i k}{n}},
	~~
	k 
	\, \in \, 
	\{ 0, \ldots , n - 1 \}.
	\non
	\ee

We are interested in designing a structured state-feedback controller, $u = - Z \psi$, that minimizes the closed-loop $\cH_2$ norm, i.e., the variance amplification from the disturbance $d$ to the regulated output $\zeta$. Since the optimal unstructured $Z$ for spatially-invariant system~\eqref{eq.spinv} is a circulant matrix~\cite{bampagdah02}, we restrict our attention to circulant feedback gains $Z$. Thus, $Z$ can also be diagonalized via a DFT and we equivalently take $x \DefinedAs \hat z$ as our optimization variable where $T^{-1} Z T = \diag \, (\hat z)$. 

For simplicity, we assume that $A$ and $Z$ are symmetric. In this case, $\hat{a}$ and $\hat z$ are real vectors and the closed-loop $\cH_2$ norm of system~\eqref{eq.spinv} takes separable form $f (x) = \sum_{k \, = \, 0}^{n - 1} f_k (x_k)$,
\[
	f_k (x_k)
	\;=\;
	\left\{
	\ba{rl}
	\dfrac{\hat q_k \, + \, \hat r_k  x_k^2}{2 \, (x_k \, - \, \hat a_k)},
	&
	x_k \,>\, \hat a_k 
	\\[0.25cm]
	\infty,
	&
	\mbox{otherwise}
	\ea
	\right.
\]
where $x_k > \hat a_k$ guarantees closed-loop stability. 
To promote sparsity of $Z$, we consider a regularized problem~\eqref{pr.split},
\be
	\ba{rl}
	\minimize\limits_{x, \, z}
	&
	f (x)
	\;+\;
	\gamma \, \| z \|_1
	\\[0.125cm]
	\subject
	&
	T x \, - \, z
	\;=\;
	0
	\ea
	\label{pr.spinv}
\ee
where $\gamma$ is a positive regularization parameter, $T$ is the inverse DFT matrix, and $z \in \bbR^n$ denotes the first row of the symmetric circulant matrix $Z$. Formulation~\eqref{pr.spinv} signifies that while it is convenient to quantify the $\cH_2$ norm in the spatial frequency domain, sparsity has to be \mbox{promoted in the physical space.} 

By solving~\eqref{pr.spinv} over a range of $\gamma$, we identify distributed controller structures which are specified by the sparsity pattern of the solutions $z^\star_\gamma$ to~\eqref{pr.spinv} at different values of $\gamma$. After selecting a controller structure associated with a particular value of $\gamma$, we solve a `polishing' or `debiasing' problem,
\be
	\ba{rl}
	\ds\minimize_{x,z}
	&
	f (x)
	\;+\;
	I_{\text{sp}(z^\star_\gamma)}(z)
	\\[0.125cm]
	\subject
	&
	T x \, - \, z
	\;=\;
	0
	\ea
	\label{pr.spinvpol}
\ee
where $I_{\text{sp}(z^\star_\gamma)}(z)$ is the indicator function associated with the sparsity pattern of $z^\star_\gamma$. The solution to this problem is the optimal controller for system~\eqref{eq.spinv} with the desired structure, i.e., the same sparsity pattern as $z^\star_\gamma$. This step is necessary because the $\ell_1$ norm in~\eqref{pr.spinv} imposes an additional penalty on $z$ that compromises closed-loop performance.

\subsubsection{Implementation}
\mbox{The elements of the gradient of $f$ are}
\[
	\dfrac{\mrd f_k (x_k)}{ \mrd x_k}
	\; = \;
	\dfrac{
	 \hat r_k x_k^2
	\, - \, 
	2\hat a_k\hat r_kx_k
	\, - \, 
	\hat q_k
	}
	 {2 \, (x_k \, - \, \hat a_k)^2}
	\]
the Hessian is a diagonal matrix with non-zero entries,
	\[
	\dfrac{\mrd^2 f_k (x_k)}{ \mrd x_k^2}
	\; = \;
	\dfrac{\hat q_k \, + \, \hat{r}_k \hat a_k^2}
	{ (x_k \, - \, \hat a_k)^3}
\]
and the proximal operator associated with the nonsmooth regularizer in~\eqref{pr.spinv} is given by soft-thresholding $\cS_{\gamma \mu}$. 

While the optimal unstructured controller can be obtained by solving $n$ uncoupled scalar quadratic equations for $x_k$, sparsity-promoting problem~\eqref{pr.spinv} is not in a separable form (because of the linear constraint) and computing the second order update~\eqref{eq.xytilde} requires solving a system of equations
\be
	\label{eq.spinveq}
	\tbt{H}{T^*}{(I - P)T}{-\mu P}
	\tbo{\tilde{x}}{\tilde y}
	\;=\;
	\tbo{\vartheta}{\theta}
\ee
Pre-multiplying by the nonsingular matrix $\blkdiag(T,I)$ and changing variables to solve for $\tilde z \DefinedAs T\tilde{x}$ brings~\eqref{eq.xytilde} into 
\[
	\tbt{T \, H \, T^{-1}}{(1/n) \, I}{I - P}{-\mu P}
	\tbo{\tilde z}{\tilde y}
	\;=\;
	\tbo{T\vartheta}{\theta}
\]
which is of the same form as~\eqref{eq.stilde}. This equation can be solved efficiently when $P$ is sparse (i.e., $\gamma \mu$ is large; cf.~Section~\ref{sec.eff}) and $\tilde{x}$ can be recovered from $\tilde z$ via FFT. 

Since the Hessian $H$ of a separable function $f$ is a diagonal matrix, the search direction can also be efficiently computed when $I - P$ is sparse (i.e., $\gamma \mu$ is small). As in Section~\ref{sec.eff}, the component of $\ty$ in the support of $P$ is determined from the bottom row of~\eqref{eq.spinveq}. The top row of~\eqref{eq.spinveq} implies $\tilde{x} = H^{-1}(\vartheta - T^* \ty)$ and substitution into the bottom row yields
\[
	(I - P) \, T \, H^{-1} \, T^* \, (I - P) \, \tilde y
	\;=\;
	\tilde{\theta}	
\]
where $\tilde{\theta} \DefinedAs (I - P) (TH^{-1}(\vartheta - T^* P \tilde y) - \theta)$ is a known vector. Thus, the component of $\tilde y$ in the support of $I - P$ can be determined by inverting a matrix whose size is determined by the support of $I - P$ and $\tilde{ x}$ is readily obtained from $\tilde y$ and $\vartheta$. The operations involving $T$ and $T^*$ can be performed via FFT; since $H$ is diagonal, multiplication by these matrices is cheap and the computational burden in solving~\eqref{eq.xytilde} again arises from a limited matrix inversion. In contrast to Section~\ref{sec.eff}, the computation of the search direction using this approach is efficient when $I - P$ is sparse, \mbox{i.e., $\gamma \mu$ is small.}

	\vsp
	
\subsubsection{Swift-Hohenberg equation}
\label{subsec.SH}
We consider the linearized Swift-Hohenberg equation~\cite{swihoh77},
\[
	\partial_t \psi(t, \xi)
	\, = \,
	(
	c I
	\,-\,
	(
	I
	\,+\,
	\partial_{\xi\xi}
	)^2)
	\,
	\psi(t,\xi)
	\,+\,
	u(t,\xi)
	\,+\,
	d(t,\xi)
\]
with periodic boundary conditions on a spatial domain $\xi \in [-\pi,\pi]$. Finite-dimensional approximation and diagonalization via the DFT (with an even number of Fourier modes $n$) yields~\eqref{eq.spinv-hat} with
$
	\hat a_k
	= 
	c - (1 - k^2)^2
$
where $k = \{ -n/2 +1, \ldots, n/2 \}$ is a spatial wavenumber. 
	
\begin{figure}[t]
\begin{center}
\vspace{-.5cm}
	\begin{tabular}{cc}
	    \!\!\!\!\!
    \!\!\!\!\!
	\subfloat[]
{
    \label{fig.sh}
	\label{}
	\begin{tabular}{rc}
    \!\!\!\!\!     \!\!\!\!\!
        \rotatebox{90}{\quad\;\footnotesize middle row of $Z$}
    \!\!\!\!\!    
    &
    \!\!\!\!\!
    \includegraphics[width= 0.45\columnwidth]{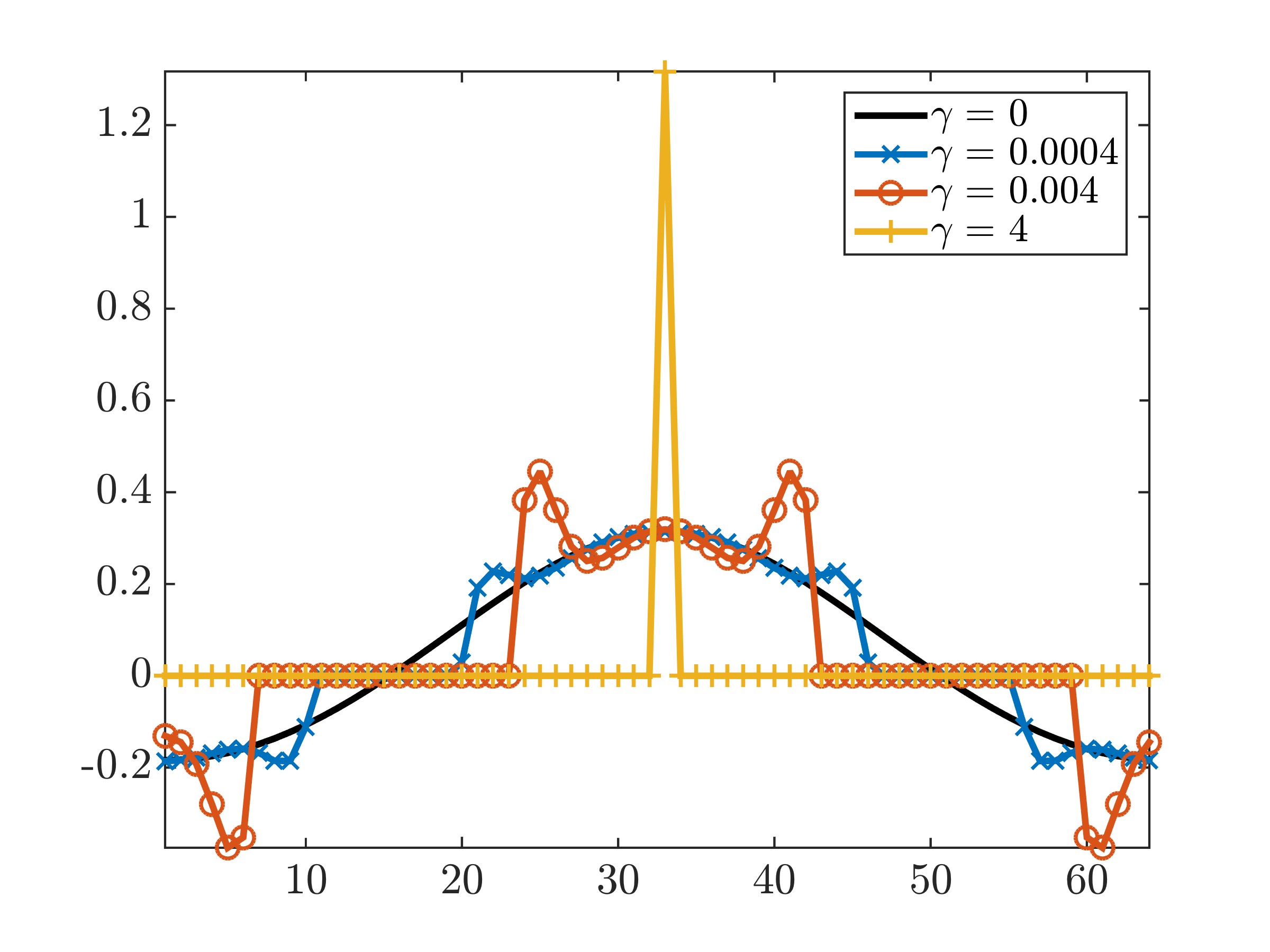}
    \!\!\!\!\!
    \!\!\!\!\!
        \\
    \!\!\!\!\!
    &
    \!\!\!\!\!
    \!\!\!\!\!
    \end{tabular}
    }
    \!\!\!\!\!
    &
        \!\!\!\!\!
  	\subfloat[]
	{
    \label{fig.gamvsp}
  \begin{tabular}{rc}
    \!\!\!\!\!
        \rotatebox{90}{\footnotesize   \quad \; sparsity level of $z^\star_\gamma$}
    \!\!\!\!\!
    &
        \!\!\!\!\!
    \!\!\!\!\!
    {	\includegraphics[width= 0.45\columnwidth]{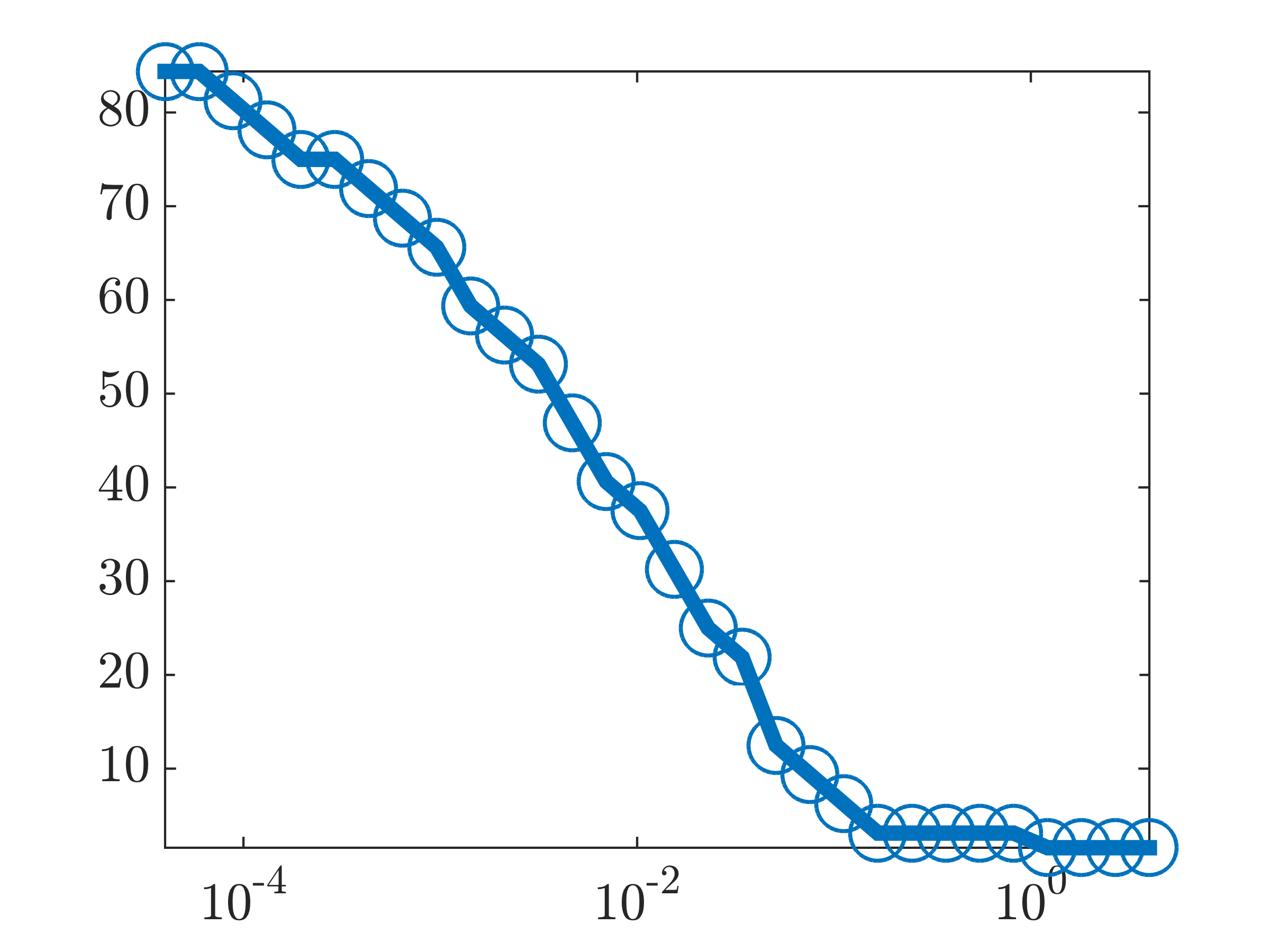}}
    \!\!\!\!\!    \!\!\!\!\!
    \\
    \!\!\!\!\!
        \!\!\!\!\!
    \!\!\!\!\!
    &
    \!\!\!\!\!
    {\footnotesize $\gamma$}
    \!\!\!\!\!
    \end{tabular}}
        \!\!\!\!\!
    \!\!\!\!\!
	\end{tabular}
\end{center}
    \caption{(a) The middle row of the circulant feedback gain matrix $Z$; and (b) the sparsity level of $z^\star_\gamma$ (relative to the sparsity level of the optimal centralized controller $z^\star_0$) resulting from the solutions to~\eqref{pr.spinv} for the linearized Swift-Hohenberg equation with $n = 64$ Fourier modes and $c = -0.01$.}
    \label{fig.trade}
\end{figure}

Figure~\ref{fig.sh} shows the optimal centralized controller and solutions to~\eqref{pr.spinv} for $c = -0.01$, $n = 64$, $Q = R = I$, and $\gamma = 4 \times 10^{-4}$, $ 4 \times 10^{-3}$, and $4$. As further illustrated in Fig.~\ref{fig.gamvsp}, the optimal solutions to~\eqref{pr.spinv} become \mbox{sparser as $\gamma$ is increased.}

\begin{figure}
\begin{center}
	\[
  \begin{tabular}{rc}   
        \rotatebox{90}{\footnotesize \!\!\!\!\!\!\!\! performance degradation (\%)}
    \!\!\!\!\!    
    &
    \!\!\!\!\!    
    {	\includegraphics[width= 0.45\columnwidth]{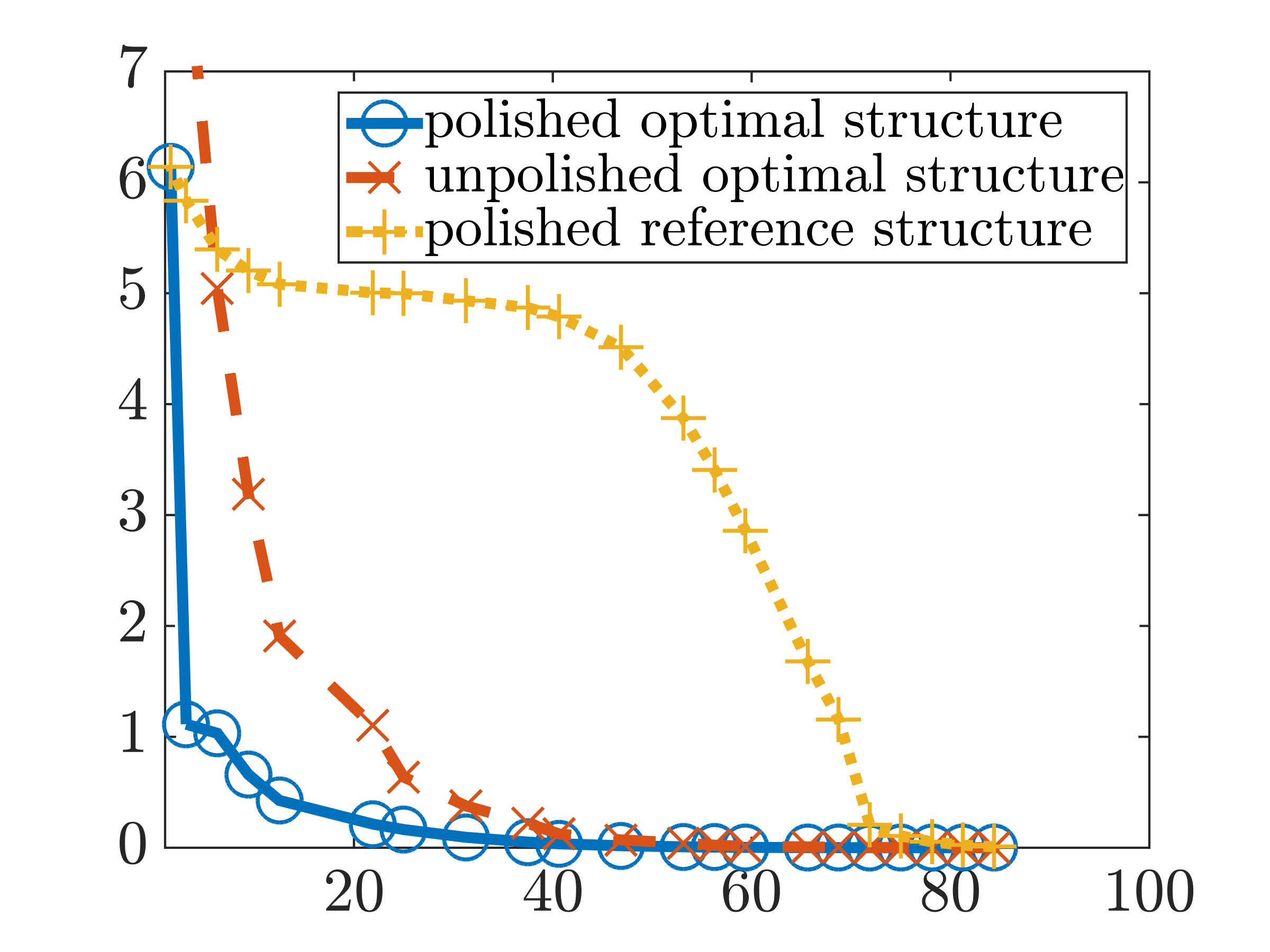}}
        \\
    \!\!\!\!\!
    &
    \!\!\!\!\!
    {\footnotesize sparsity level of $z$ (\%)}
    \end{tabular}
    \]
    \end{center}
    \caption{Performance degradation (in percents) of structured controllers relative to the optimal centralized controller: polished optimal structured controller obtained by solving~\eqref{pr.spinv} and~\eqref{pr.spinvpol} (\tc{matlabblue}{---$\circ$---}); unpolished optimal structured controller obtained by solving only~\eqref{pr.spinv} (\tc{matlabred}{- -$\times$- -}); and optimal structured controller obtained by solving~\eqref{pr.spinvpol} for an {\em a priori\/} specified nearest neighbor reference structure (\tc{matlaborange}{$\cdots+\cdots$}).
    }
    \label{fig.spvf}
\end{figure}

In Fig.~\ref{fig.spvf}, we demonstrate the utility of using regularized problems to navigate the tradeoff between controller performance and structure, an approach pioneered by~\cite{linfarjovTAC13admm}. The polished optimal structured controllers (\tc{matlabblue}{---$\circ$---}) were designed by first solving~\eqref{pr.spinv} to identify an optimal structure and then solving~\eqref{pr.spinvpol} to further improve the closed-loop performance. To illustrate the importance of polishing step~\eqref{pr.spinvpol}, we also show the closed-loop performance of unpolished optimal structured controllers (\tc{matlabred}{- -$\times$- -}) resulting from~\eqref{pr.spinv}. Finally, to evaluate the controller {\em structures\/} identified by~\eqref{pr.spinv}, we show the closed-loop performance of polished `reference' structured controllers (\tc{matlaborange}{$\cdots + \cdots$}). Instead of solving~\eqref{pr.spinv}, the reference structures are {\em a priori\/} specified as nearest neighbor symmetric controllers of the same cardinality as controllers resulting from~\eqref{pr.spinv}. Among controllers with the same number of nonzero entries, the polished optimal structured controller consistently achieves the best closed-loop performance.

We compare the computational efficiency of our approach with the proximal Newton method~\cite{leesunsau14} and ADMM~\cite{wujovSCL17}. The proximal Newton method requires solving a LASSO subproblem~\eqref{eq.proxnewt}, for which we employ SpaRSA~\cite{wrinowfig09}. \tcbl{ADMM was implemented using the adaptive stepsize selection discussed in~\cite[Section 3.4.1]{boyparchupeleck11} from~\cite{hexyanwan00,wanlia01}.}
Since $A$ is circulant, the $x$-minimization step in ADMM~\eqref{eq.ADMM} requires solving $n$ uncoupled cubic scalar equations. In general, when $A$ is block circulant, the DFT only block diagonalizes the dynamics and thus the $x$-minimization step has to be solved via an iterative procedure~\cite{wujovSCL17}.

\begin{figure}
\begin{center}
\[
	\begin{tabular}{rc}
    \!\!\!\!\!
        \rotatebox{90}{\footnotesize \quad\quad\quad\quad time}
    \!\!\!\!\!
    &
    \!\!\!\!\!
    {	\includegraphics[width= 0.45\columnwidth]{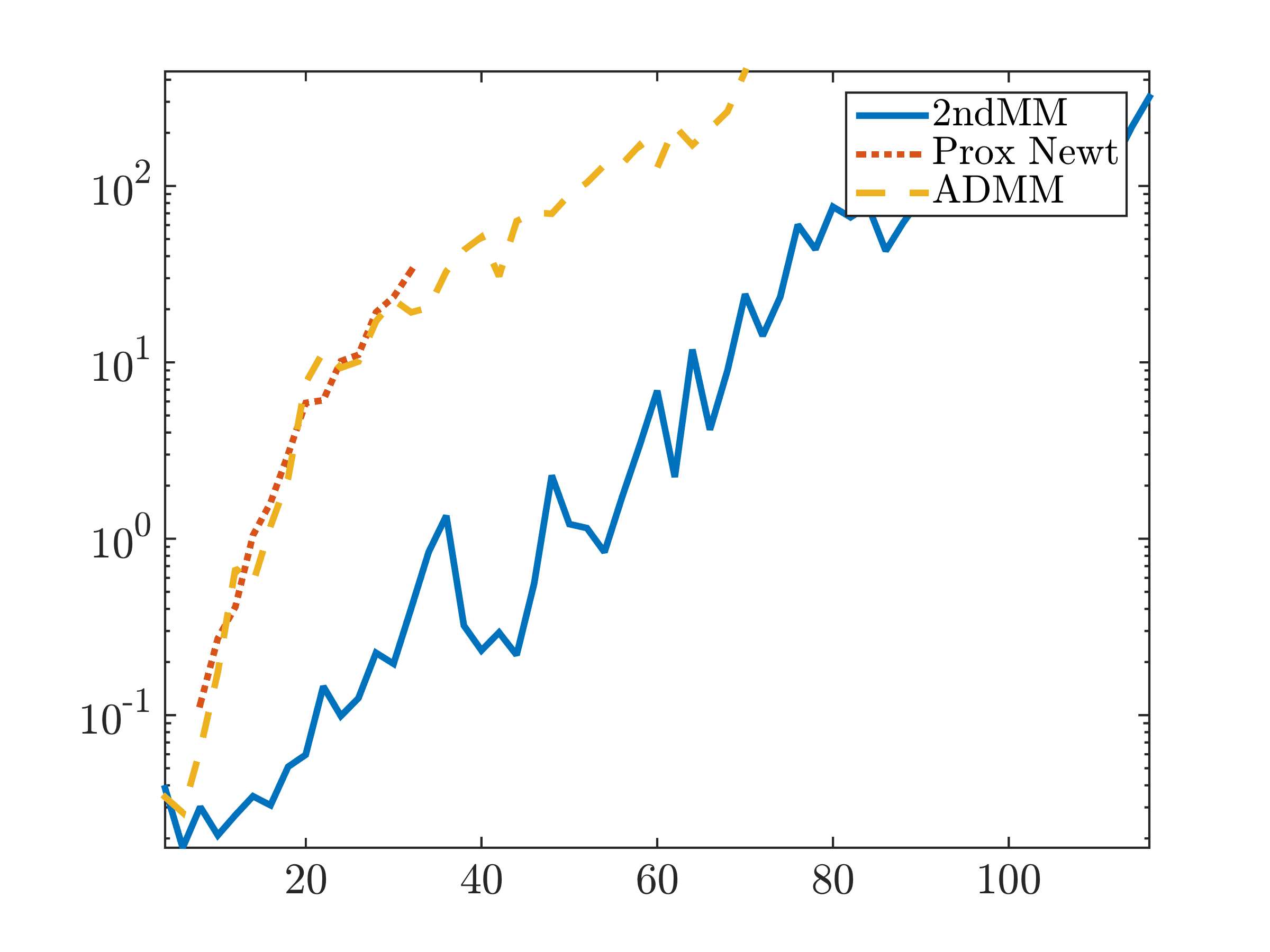}}
        \\
    \!\!\!\!\!
    \!\!\!\!\!
    &
    \!\!\!\!\!
    {\footnotesize $n$}
    \!\!\!\!\!
    \end{tabular}
\]
\end{center}
    \caption{Total time to compute a solution to~\eqref{pr.spinv} with $\gamma = 0.004$ using our algorithm (2ndMM), proximal Newton, and ADMM.}
    \label{fig.sht}
\end{figure}

\begin{figure}[t]
\begin{center}
\vspace{-.5cm}
\begin{tabular}{cc}
    \!\!\!\!\!
    \!\!\!\!\!
\subfloat[][]
{
\label{fig.shinner}
  \begin{tabular}{rc}
    \!\!\!\!\!    \!\!\!\!\!
        \rotatebox{90}{\footnotesize   ~~\;\quad avg iteration time}
    \!\!\!\!\!  
    &
        \!\!\!\!\!
    {	\includegraphics[width= 0.45\columnwidth]{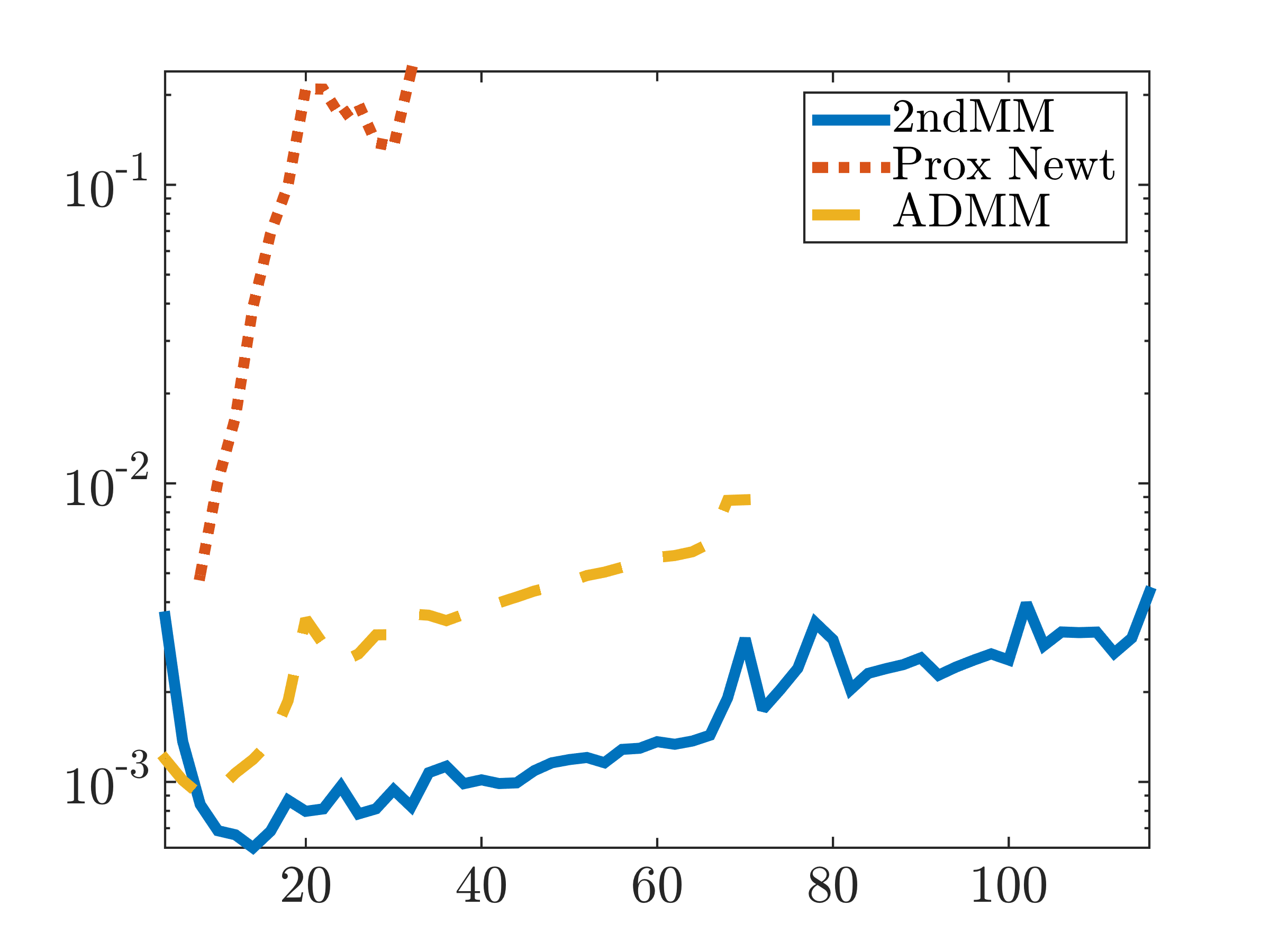}}
    \!\!\!\!\!    \!\!\!\!\!
    \\
    \!\!\!\!\!
    &
    \!\!\!\!\!
    {\footnotesize $n$}
    \!\!\!\!\!
    \end{tabular}
}
    \!\!\!\!\!
&
    \!\!\!\!\!
    \subfloat[][]
    {
    \label{fig.shiter}
  \begin{tabular}{rc}
    \!\!\!\!\!    
        \rotatebox{90}{\footnotesize   \quad\; number of iterations}
    \!\!\!\!\!    
    &
    \!\!\!\!\!    \!\!\!\!\!
    {	\includegraphics[width= 0.45\columnwidth]{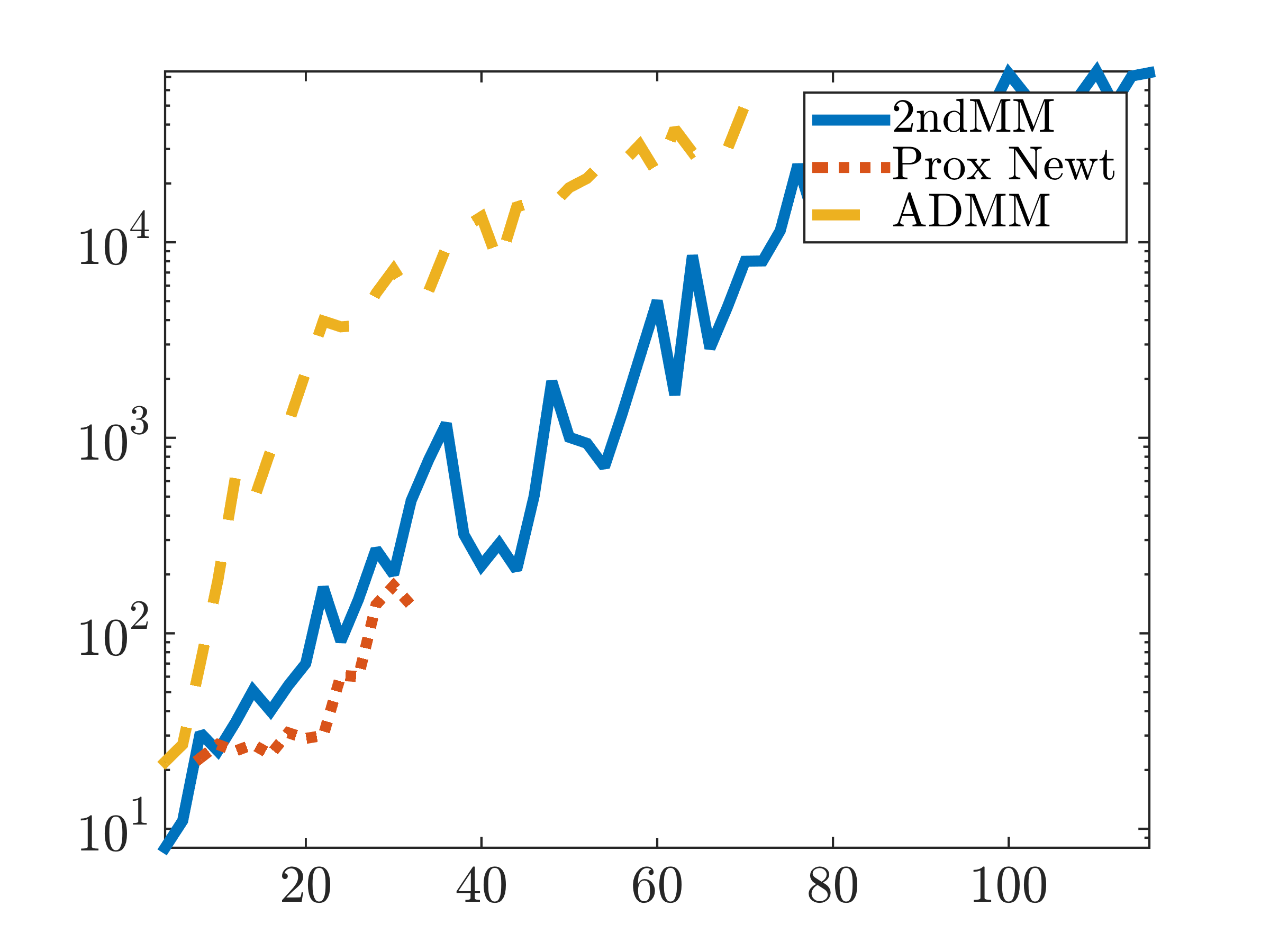}}
        \\
    \!\!\!\!\!
    \!\!\!\!\!
    &
    \!\!\!\!\!
    {\footnotesize $n$}
    \!\!\!\!\!
    \end{tabular}
    \!\!\!\!\!
    \!\!\!\!\!
}
\end{tabular}
\end{center}
    \caption{Comparison of (a) times to compute an iteration (averaged over all iterations); and (b) numbers of iterations required to solve~\eqref{pr.spinv} with $\gamma = 0.004$.
    }
        \!\!\!\!\!
    \!\!\!\!\!
    \label{fig.shtdet}
\end{figure}

Figure~\ref{fig.sht} shows the time to solve~\eqref{pr.spinv} with $\gamma = 0.004$ using our method, proximal Newton, and ADMM. Our algorithm and ADMM were stopped when the primal and dual residuals were below $1 \times 10^{-8}$. The proximal Newton method was stopped when the norm of the difference between two consecutive iterates was smaller than $1 \times 10^{-8}$.  In Fig.~\ref{fig.shtdet}, we show the per iteration cost and the total number of iterations required to find the optimal solution using each method. Our algorithm clearly outperforms proximal Newton and ADMM.
	
 Although proximal Newton requires a similar number of iterations, the LASSO subproblem~\eqref{eq.proxnewt} that determines its search direction is much more expensive; this increases the computation cost of each iteration and slows the overall algorithm. Moreover, for larger problem sizes, the proximal Newton method struggles with finding a stabilizing search direction because $\nabla^2 f$ seems to bring it away from the set of stabilizing feedback gains. It appears that our method circumvents this issue because its iterates lie in a larger lifted space in which stability is easier to enforce via backtracking.

On the other hand, while the $x$- and $z$- minimization steps in ADMM are quite efficient, as a first order method ADMM requires a large number of iterations to reach high-accuracy solutions. Our algorithm achieves better performance because its use of second order information leads to relatively few iterations and the structured matrix inversion leads to efficient computation of the search direction.

	\vspace*{-1ex}
\section{Connections and discussion}
\label{sec.disc}

The proximal augmented Lagrangian $\cL_\mu(x;y)$ is obtained by constraining $\cL_\mu(x,z;y)$ to the manifold
\[
	{\cal Z}
	\; \DefinedAs \;
	\{
	(x,z^{\star}_\mu;y)
	~|~
	z^\star_\mu
	\;=\;
	\argmin_{z}
	\,
	\cL_\mu(x,z;y)
	\}
	\; = \;
	\{
	(x,z^{\star}_\mu;y)
	~|~
	Tx \, + \, \mu y
	\, \in \,
	z^\star_\mu
	\, + \,
	\mu
	\, 
	\partial_C g (z^\star_\mu)
	\}
\]
which results from the explicit minimization over the auxiliary variable $z$. Herein, we interpret the second order search direction as a linearized update to the KKT conditions for problem~\eqref{pr.split} and discuss connections to \mbox{the alternative algorithms.}

	\vspace*{-2ex}
\subsection{Second order updates as linearized KKT corrections}

The second order update~\eqref{eq.xytilde} can be viewed as a first order correction to the KKT conditions for optimization problem~\eqref{pr.split},
\be
	\ba{rcl}
	0
	&\!\!=\!\!&
	\nabla f(x) 
	\,+\,
	T^Ty
	\\[0.cm]
		0
	&\!\!=\!\!&
	Tx
	\,-\,
	z
	\\[0.cm]
	0
	&\!\! \in \!\!&
	\partial_C g(z)
	\,-\,
	y.
	\ea
	\label{eq.KKT}
\ee
Substitution of $z^\star_\mu$ into~\eqref{eq.KKT} makes the last two conditions redundant; when combined with the definition of the manifold $\cal Z$, $Tx = z$ implies $y \in \partial_C g(z)$ and $y \in \partial_C g(z)$ implies $Tx = z$.

Multiplying equation~\eqref{eq.xytilde} for the second order update with the nonsingular matrix~\eqref{eq.Pi} recovers the equivalent expression
	\be
	\tbt{\!\!\Hes\!\!}{\!\!T^T\!\!}{\!\!(I \,-\, P)T\!\!}{\!\!-\mu P\!\!}
	\tbo{\tx^k}{\ty^k}
	= 
	-
	\tbo{\gr(x^k) + T^Ty^k}{r^k}
	\label{eq.wtilde}
	\ee
where $r^k \DefinedAs Tx^k - z^\star_\mu(x^k;y^k) = Tx^k - \prox_{\mu g} (Tx^k + \mu y^k)$ is the primal residual of~\eqref{pr.split}. Thus,~\eqref{eq.wtilde} describes a first order correction to the first and second KKT conditions in~\eqref{eq.KKT}.

	\vspace*{-2ex}
\subsection{Connections with other methods}

We now discuss broader implications of our framework and draw connections to the existing methods for solving versions of~\eqref{pr}. Many techniques for solving composite minimization problems of the form~\eqref{pr} can be expressed in terms of functions embedded in the augmented Lagrangian; see Table~\ref{tab.embed}. Trivially, the  objective function in~\eqref{pr} corresponds to $\cL_\mu(x,z;y)$ over the manifold $z = Tx$. The Lagrange dual of a problem equivalent to~\eqref{pr.split},
\be
	\label{pr.str}
	\ba{rl}
	\minimize\limits_{x, \, z}
	&
	f(x) \;+\; g(z) \;+\; \tfrac{1}{2\mu} \, \norm{Tx \,-\, z}^2
	\\
	\subject
	&
	Tx \,-\, z
	\;=\;
	0
	\ea
\ee
is recovered by collapsing $\cL_\mu(x,z;y)$ onto the intersection of the $z$-minimization manifold $\cal Z$ and the $x$-minimization manifold,
\[
	{\cal X}
	\; \DefinedAs \;
	\{
	(x^\star,z;y)
	~|~
	x^\star 
	\, = \,
	\argmin_{x}
	\,
	\cL_\mu(x,z;y)
	\}
	\; = \;
	\{
	(x^\star,z;y)
	~|~
	\mu
	\,
	(
	\nabla f (x^\star)
	\, + \, 
	T^T y
	)
	\, = \,
	T^T(z \,-\, Tx^\star)
	\}.
\]
The Lagrange dual of~\eqref{pr.split} is recovered from the Lagrange dual of~\eqref{pr.str} in the limit $\mu \to \infty$.

	\vsp
	
\subsubsection{MM and ADMM}

MM implements gradient ascent on the dual of~\eqref{pr.str} by collapsing $\cL_\mu(x,z;y)$ computationally onto $\cal X \cap \cal Z$. The joint $(x,z)$-minimization step in~\eqref{eq.MM} evaluates the Lagrange dual at discrete iterates $y^k$ by finding the corresponding $(x,z)$-pair on $\cal X \cap \cal Z$; i.e., the iterate $(x^{k+1},z^{k+1};y^k)$ generated by~\eqref{eq.MM} lies on the manifold $\cal X \cap \cal Z$. Note that, in this form, joint $(x,z)$-minimization is a challenging nondifferentiable optimization problem.

ADMM avoids this challenge by collapsing $\cL_\mu(x,z;y)$ onto $\cal X$ and $\cal Z$ separately. While the underlying $x$- and $z$-minimization subproblems in ADMM are relatively simple, the iterate $(x^{k+1},z^{k+1}; y^k)$ generated by~\eqref{eq.ADMM} does not typically lie on the $\cal X \cap \cal Z$ manifold. Thus ADMM does not represent gradient ascent on the dual of~\eqref{pr.str}, causing looser theoretical guarantees and often worse practical performance. \tcbl{Although second order methods may be employed within ADMM (e.g., in \eqref{eq.ADMM} to solve the {\em inner\/} $x$-minimization problem), ADMM is still a first order algorithm as the outer iteration consists of gradient updates to the dual variable.}

\begin{table*}
\begin{center}
  \begin{tabular}{ | c | c | c | c | c |}
    \hline
    {\bf$x$} & {\bf$z$} & {\bf$y$} & {\bf function} & {\bf methods} \\ \hline \hline
    \multicolumn{2}{|c|}{$z = Tx$} & - & objective function of~\eqref{pr} & subgradient iteration~\cite{ber99} \\ \hline
    - & - & - & augmented Lagrangian & ADMM~\cite{linfarjovTAC13admm,boyparchupeleck11} \\ \hline
    $\cal X$ & $\cal Z$ & - & Lagrange dual of~\eqref{pr.str} & Dual ascent (if explicit expression available) \\ \hline
    - & $\cal Z$ & - & proximal augmented Lagrangian & MM~\cite{dhikhojovTAC19}, Arrow-Hurwicz-Uzawa~\cite{dhikhojovTAC19}, {\bf Second order primal-dual} \\ \hline
    - & $\cal Z$ & $\cal Y$ & Forward-Backward Envelope & Forward-Backward Truncated Newton~\cite{patstebem14}, related to prox. gradient~\cite{patstebem14} \\ \hline
      \end{tabular}
\end{center}
\caption{Summary of different functions embedded in the augmented Lagrangian of~\eqref{pr.split} and methods for solving~\eqref{pr} based on these functions. \tcbl{Each function corresponds to the augmented Lagrangian on a different manifold. The columns labeled $x$, $y$, and $z$ contain the manifolds to which the corresponding variables in the augmented Lagrangian must be constrained in order to obtain the function in that row.}
}
\label{tab.embed}
\end{table*}

By collapsing the augmented Lagrangian onto $\cal Z$, the proximal augmented Lagrangian~\eqref{eq.alprox} allows us to express the $(x,z)$-minimization step in MM as a tractable continuously differentiable problem~(cf.\ Theorem~\ref{thm.diff}). This avoids challenges associated with ADMM and it does not increase computational complexity in the $x$-minimization step in~\eqref{eq.MM} relative to ADMM when using first order methods. We finally note that unlike the Rockafellar's {\em proximal method of multipliers\/}~\cite{roc76} which applies the proximal point algorithm to the primal-dual optimality conditions, our framework {\em reformulates\/} the {\em standard method of multipliers\/} and develops second order algorithm to solve nonsmooth composite optimization problems.

\subsubsection{Second order methods}

We identify the saddle point of $\cL_\mu(x,z;y)$ by forming second order updates to $x$ and $y$ along the $\cal Z$ manifold. Just as Newton's method approximates an objective function with a convex quadratic function, we approximate $\cL_{\mu}(x;y)$ with a quadratic saddle function.

Constraining the dual variable itself yields connections with other methods. When $T = I$,~\eqref{eq.KKT} implies that the optimal dual variable is given by $y^\star = -\nabla f(x^\star)$, so it is natural to collapse the augmented Lagrangian onto the manifold
	$
	{\cal Y}
	\DefinedAs
	\{
	(x,z;y^\star)
	~|~
	y^\star 
	= 
	-
	\nabla f(x)
	\}.
	$
The augmented Lagrangian over the manifold $\cal Z \cap \cal Y$ corresponds to the Forward-Backward Envelope (FBE) introduced in~\cite{patstebem14}. The proximal gradient algorithm with step size $\mu$ can be recovered from a {variable metric} gradient iteration on the FBE~\cite{patstebem14}. In~\cite{patstebem14,stethepat17,thestepat18,stethepat18}, the approximate line search and quasi-Newton methods based on the FBE were developed to solve~\eqref{pr}. Since the Hessian of the FBE involves third order derivatives of $f$, these techniques employ either truncated- or quasi-Newton methods. 

	The approach advanced in the current paper applies a second order method to the augmented Lagrangian that is constrained over the larger manifold $\cal Z$. Relative to alternatives, our framework offers several advantages. First, while the FBE is in general nonconvex function of $x$, $\cL_\mu(x;y)$ is always convex in $x$ and concave in $y$. Furthermore, we can compute the Hessian exactly using only second order derivatives of $f$ and its structure allows for efficient computation of the search direction. Finally, our formulation allows us to leverage recent advances in second order methods for augmented Lagrangian methods, e.g.,~\cite{gilrobCOA12,armbenomhpat14,armomh15}.

	\vspace*{-2ex}
\section{Concluding remarks}
	\label{sec.conc}
	
We have developed a second order primal-dual method for nonsmooth convex composite optimization. We establish global exponential convergence of the corresponding differential inclusion in continuous-time and use the primal-dual augmented Lagrangian as a merit function in a discrete-time implementation. Our globally convergent customized algorithm exhibits asymptotic (quadratic) superlinear convergence rate when the proximal operator associated with nonsmooth regularizer is (strongly) semismooth. We use the $\ell_1$-regularized least squares and spatially-invariant control problems to demonstrate competitive performance of our algorithm relative to the available state-of-the-art alternatives.

Future research will focus on \tcbl{weakly convex problems,} nonconvex problems, the development of inexact second order methods, and application of our framework to problems in which the codomain of the regularizer has a larger dimensions than the domain of the optimization variable. In particular, inexact but structured approximations of the Hessian can lead to efficient methods that may even be convenient for distributed implementation.

	\vspace*{-2ex}
\appendix

\subsection{Algorithm based on $V(w)$ as a merit function} \label{app.alg}

Since Theorem~\ref{thm.inc} establishes global convergence of the differential inclusion, one algorithmic approach is to implement a Forward Euler discretization of differential inclusion~\eqref{eq.diffinc}. A natural choice of merit function is the Lyapunov function $V$ defined in~\eqref{eq.V}. A simple corollary of Theorem~\ref{thm.asym} shows that the second order update~\eqref{eq.xytilde} is a descent direction for $V$.

	\vsp
	
\begin{corollary} \label{cor.desc}
	The second order update~\eqref{eq.xytilde} is a descent direction for the merit function $V$ defined in~\eqref{eq.V}.
\end{corollary}

	\vsp
	
\begin{IEEEproof}
	Follows from~\eqref{eq.gexp2} in Theorem~\ref{thm.asym}.
\end{IEEEproof}

	\vsp

Corollary~\ref{cor.desc} enables the use of a backtracking Armijo rule for step size selection. A natural choice of stopping criterion for such an algorithm is a condition on the size of the primal and dual residuals. Moreover, Proposition~\ref{prop.quad} suggests fast asymptotic convergence when $\prox_{\mu g}$ is semismooth. The example in Fig.~\ref{fig.sommv} verifies this intuition when solving LASSO to a threshold of $1 \times 10^{-8}$ for the primal and dual residuals.

\begin{algorithm}
\caption{Second order primal-dual algorithm for nonsmooth composite optimization based on discretizing~\eqref{eq.diffinc}.}
\label{alg.v}
\begin{algorithmic}
\STATE \textbf{input:} Initial point $x_0$, $y_0$, backtracking constant $\alpha \in (0,1)$, Armijo parameter $\sigma \in (0,1)$, and stopping tolerances $\varepsilon_1$, $\varepsilon_2$.
 \vspace*{-0.15cm}
\\
\textbf{While:} $\norm{Tx^k - \prox_{\mu g}(Tx^k + \mu y^k)} > \varepsilon_1$ or

\vspace*{0.05cm}

\noindent\phantom{\textbf{While:}} $\norm{\nabla f (x^k) \,-\, T^Ty^k} > \varepsilon_2$

\vspace{0.1cm}
~\,\quad \textbf{Step $1$:} Compute $\tw^k$ as defined in~\eqref{eq.xytilde}

~\,\quad \textbf{Step $2$:} Choose the smallest $j \in \mathbb{Z}_{+}$ such that
\[
	V(w^k + \alpha^j \tw^k)
	\;\leq\;
		V(w^k)
	\,-\,
	\sigma
	\,
	\alpha^j
	\,
	\norm{\nabla\cL_\mu(w^k)}^2
\]
~\,\quad \textbf{Step $3$:} Update the primal and dual variables
\[
	w^{k+1}
	\, = \,
	w^k
	\, + \,
	\alpha^j \tw^k
\]
\end{algorithmic}
\end{algorithm}

However, such an implementation would require a fixed penalty parameter $\mu$, which typically has a large effect on the convergence speed of augmented Lagrangian algorithms and is difficult to select {\em a priori\/}.  Moreover, stability of the solution to a differential equation does not always imply stability of its discretization.

\subsubsection*{An example}
\begin{figure}
\begin{center}
\[
	\begin{tabular}{rc}
    \!\!\!\!\!
        \rotatebox{90}{\footnotesize \quad\quad\quad  $\norm{x - x^\star}$}
    \!\!\!\!\!
    &
    \!\!\!\!\!
    {	\includegraphics[width= 0.45\columnwidth]{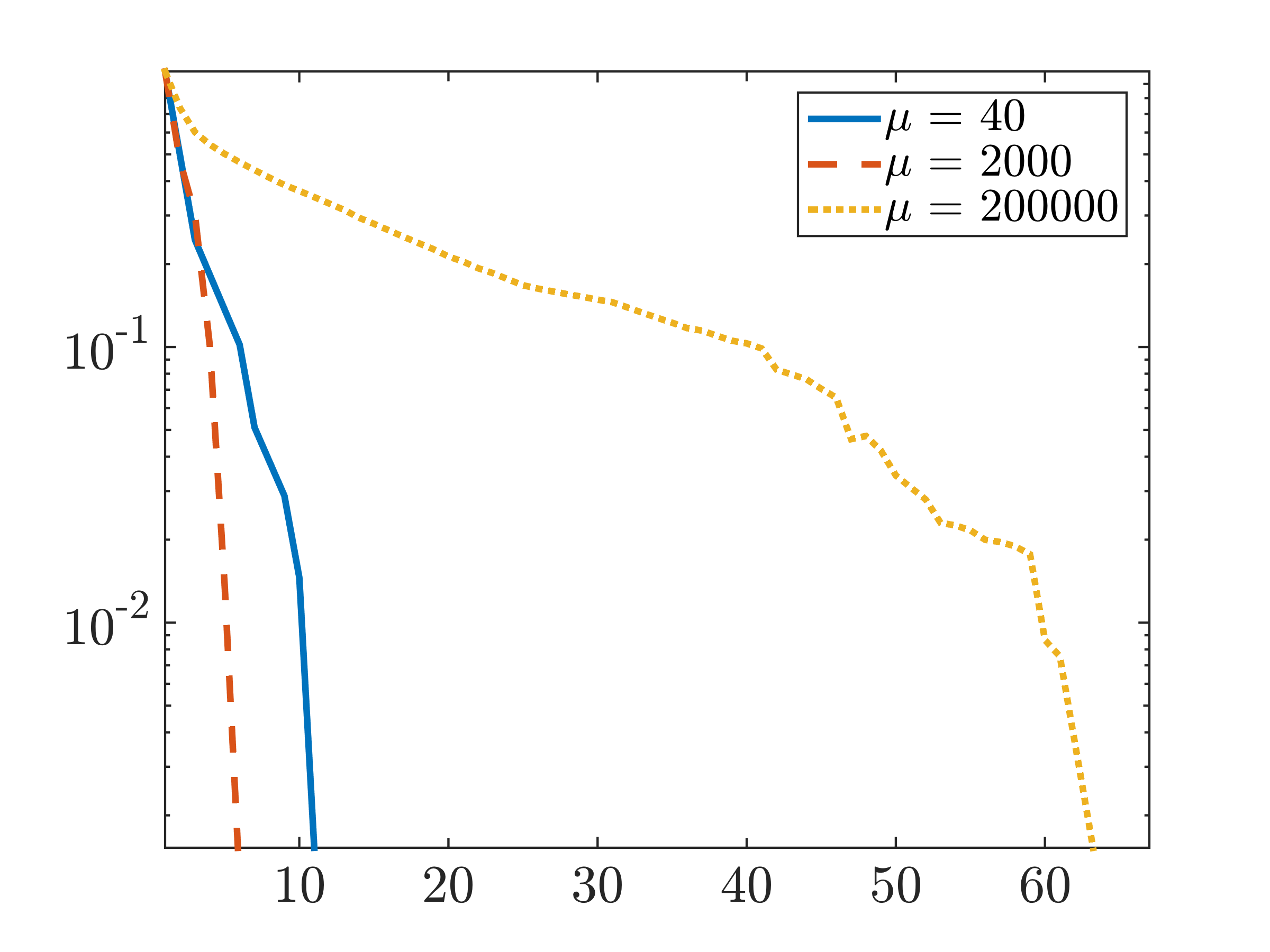}}
        \\
    \!\!\!\!\!
    \!\!\!\!\!
    &
    \!\!\!\!\!
    {\footnotesize Iteration}
    \end{tabular}
\]
\end{center}
    \caption{Distance from the optimal solution as a function of iteration number when solving LASSO using Algorithm~\ref{alg.v} for different values of $\mu$.}
    \label{fig.sommv}
\end{figure}

We implement Algorithm~\ref{alg.v} to solve the LASSO problem~\eqref{eq.lasso} studied in Section~\ref{sec.ex}. The problem instance was randomly generated with $n = 500$, $m = 1000$, and $\gamma = 0.85\gamma_{\max}$. Figure~\ref{fig.sommv} illustrates the quadratic asymptotic convergence of Algorithm~\ref{alg.v} and a strong influence of $\mu$.

	\vspace*{-2ex}
\subsection{Bounding $\nabla f(x)$ for Algorithm~\ref{alg.sqp}}\label{sec.appf}

A bounded set $\cC_f$ containing the solution to~\eqref{pr} can always be identified from any point $\bar x$, $\nabla f(\bar x)$,  any element of the subgradient $\partial g(\bar x)$, and a lower bound on the parameter of strong convexity.
\begin{lemma} \label{lem.bdset}
	Let Assumption~\ref{ass.func} hold and let $f$ be $m_f$-strongly convex.
	Then, for any $\bar x$, the optimal solution to~\eqref{pr} lies within the ball of radius $\tfrac{2}{m_f}\norm{\nabla f(\bar x) \,+\, T^T\partial g(T\bar x)}$ centered at $\bar x$.
\end{lemma}

\begin{IEEEproof}
	Given any point $\bar x$, strong convexity of $f$ and convexity of $g$ imply that,
	\vspace{-0.2cm}
	\begin{equation}
		\label{eq.ineq}
		f(x) \, + \, g(Tx)
		\,-\,
		f(\bar x) \, - \, g(T\bar x) 
		\;\geq\;
		\inner{\nabla f(\bar x) \,+\, T^T\partial g(T\bar x)}{x \,-\, \bar x}
		\, +\,
		\tfrac{m_f}{2}
		\norm{x \,-\, \bar x}^2
	\end{equation}
	for all $x$ where $\partial g(T\bar x)$ is any member of the subgradient of $g(T\bar x)$. For any $x$ with
	$
	\norm{x - \bar x} \geq \tfrac{2}{m_f}\norm{\nabla f(\bar x) + T^T\partial g(T\bar x)},
	$
	 the right-hand side of~\eqref{eq.ineq} must be nonnegative which implies that $f(x) + g(Tx) \geq f(\bar x) + g(T\bar x)$ and thus that $x$ cannot solve~\eqref{pr}.
	\end{IEEEproof}
	
For any strongly convex function $f$ and convex function $g$, a function $\hat f$ can be selected such that
\[
	\argmin\, f(x) \,+\, g(Tx)
	\;=\;
	\argmin\, \hat f(x) \,+\, g(Tx).
\]
Here, $\hat f$ is identical to $f$ in some closed set $\cC_f$ containing $x^\star$,
	$
	\hat f(x)
	\DefinedAs
	\{
	f(x),
	\,
	x \in \cC_f;
	$
	$
	\tilde f(x),
	\,
	x \not\in \cC_f
	\},
	$
and $\tilde f(x)$ is chosen such that $\hat f(x)$ is convex and twice differentiable, $\nabla \hat f(x)$ is uniformly bounded, and $\nabla^2 \hat f(x) \succ 0$. Lemma~\ref{lem.bdset} implies that a set $\cC_f$ that contains the optimal solution $x^\star$ of~\eqref{pr} can be identified {\em a priori\/} from any given point $\bar x$. Since $\hat f(x)$ satisfies all the conditions of Theorem~\ref{thm.conv}, Algorithm~\ref{alg.sqp} can be used to solve~\eqref{pr} and, since $\cC_f$ contains $x^\star$, its optimal solution also solves~\eqref{pr}.

\subsubsection*{An example}  When $x \in \bbR$, $f(x) = \tfrac{1}{2}x^2$, and $\cC_f = [-1, 1]$, a potential choice of $\tilde f(x)$ is given by,
	$
	\tilde f(x)
	=
	\left\{
	-2x + \mre^{x + 1} - 2.5,
	~
	x \leq -1; 
	\right.
	$
	$
	\left.
	2x + \mre^{-x + 1}  - 2.5,
	~
	x \geq 1
	\right\}.
	$
For this choice of $\tilde f$, the gradient of $\hat f$ is continuous and bounded,
	$
	\nabla \hat f(x)
	=
	\{
	-2 + \mre^{x + 1},
	\,
	x \leq -1;
	x,
	\,
	x \in [-1,1];
	2 - \mre^{-x + 1},
	\,
	x \geq 1
	\}
	$
and the Hessian of $\hat f$ is determined by
	$
	\nabla^2 \hat f(x)
	=
	\{
	\mre^{x + 1},
	\,
	x \leq -1;
	1,
	\,
	x \in [-1,1];
	\mre^{-x + 1},
	\,
	x \geq 1
	\}.
	$

\end{document}